\documentclass[12pt]{article}

\usepackage{amssymb}
\usepackage{amsmath}
\usepackage{cite}
\usepackage{amsthm}

\begin{document}
\renewcommand{\citeleft}{{\rm [}}
\renewcommand{\citeright}{{\rm ]}}
\renewcommand{\citepunct}{{\rm,\ }}
\renewcommand{\citemid}{{\rm,\ }}

\newcounter{abschnitt}
\newtheorem{satz}{Theorem}
\newtheorem{theorem}{Theorem}[abschnitt]
\newtheorem{koro}[theorem]{Corollary}
\newtheorem{prop}[theorem]{Proposition}
\newtheorem{lem}[theorem]{Lemma}

\newcounter{saveeqn}
\newcommand{\alpheqn}{\setcounter{saveeqn}{\value{abschnitt}}
\renewcommand{\theequation}{\mbox{\arabic{saveeqn}.\arabic{equation}}}}
\newcommand{\reseteqn}{\setcounter{equation}{0}
\renewcommand{\theequation}{\arabic{equation}}}

\hyphenation{convex} \hyphenation{bodies}

\sloppy

\phantom{a}

\vspace{-1.1cm}

\begin{center}

\begin{large} {\bf Harmonic Analysis of Translation Invariant Valuations} \\[0.7cm] \end{large}

\textsc{Semyon Alesker, Andreas Bernig and Franz E. Schuster}
\end{center}

\vspace{-1.1cm}

\begin{quote}
\footnotesize{ \vskip 1truecm\noindent {\bf Abstract.} The
decomposition of the space of continuous and translation
invariant valuations into a sum of $\mathrm{SO}(n)$ irreducible
subspaces is obtained. A reformulation of this result in terms of
a Hadwiger type theorem for \linebreak continuous translation
invariant and $\mathrm{SO}(n)$-equivariant tensor valuations is
also given. As an application, symmetry properties of rigid
motion invariant and homogeneous bivaluations are established and
then used to prove new inequalities of Brunn--Minkowski type for
convex body valued valuations. }
\end{quote}

\vspace{0.6cm}

\centerline{\large{\bf{ \setcounter{abschnitt}{1}
\arabic{abschnitt}. Introduction and statement of main results}}}

\vspace{0.6cm} \reseteqn \alpheqn \setcounter{theorem}{0}

Let $V$ be an $n$-dimensional Euclidean vector space and let $A$
be an abelian semigroup. A function $\phi$ defined on convex
bodies (compact convex sets) in $V$ and taking values in $A$ is
called a {\it valuation}, or {\it additive}, if
\[\phi(K) + \phi(L) = \phi(K \cup L) + \phi(K \cap L)\]
whenever $K, L$ and $K \cup L$ are convex.

The most important cases are $A=\mathbb{R}$ or $\mathbb{C}$
(scalar valued valuations), $A=\mathrm{Sym}^k V$ (tensor
valuations) and $A=\mathcal{K}^n$, the semigroup of convex bodies
in $V$ with the Minkowski addition (Minkowski valuations).

Scalar valued valuations play an important role in integral
geometry. \linebreak Hadwiger characterized in
\textbf{\cite{hadwiger51}} the continuous Euclidean motion
invariant valuations. Almost all classical integral-geometric
formulas can be reduced to this landmark result. For
generalizations of this idea in different directions, we refer to
\textbf{\cite{Klain:Rota, Alesker99, Alesker01, Bernig09,
bernigfu10, fu06, centro, McMullen93}}.

Tensor valuations were studied by McMullen \textbf{\cite{McMullen97}}, the first author
\textbf{\cite{Alesker97}} and Ludwig
\textbf{\cite{Ludwig:matrix}}. Recently, a full set of kinematic formulas for tensor valuations was obtained
by Hug, Schneider and R. Schuster \textbf{\cite{hug_schneider_schuster_a, hug_schneider_schuster_b}}.

The best known example of a Minkowski valuation is the projection
body. This central notion from affine geometry has many
applications in several areas such as geometric tomography,
stereology, computational geometry, optimization or functional
analysis. For a systematic study of Minkowski valuations, we
refer to \textbf{\cite{ludwig02, Ludwig:Minkowski, haberl08,
haberl09, kiderlen05, Ludwig:matrix, Ludwig06, Schu09}} and the
references therein.

In this article, we contribute to these three different directions
in the \linebreak theory of valuations. Our main result may be
stated in the language of scalar valued valuations or in the
language of tensor valuations. For simplicity, we assume
throughout this paper that $n \geq 3$, even if most of the
results also hold true for $n=2$.

A valuation $\phi$ is called {\it translation invariant} if
$\phi(K+x)=\phi(K)$ for all $x \in V$ and $K \in \mathcal{K}^n$
and $\phi$ is said to have {\it degree} $i$ if $\phi(t
K)=t^i\phi(K)$ for all $K \in \mathcal{K}^n$ and $t > 0$. We call
$\phi$ {\it even} if $\phi(-K)=\phi(K)$ and {\it odd} if
$\phi(-K)=-\phi(K)$ for all $K \in \mathcal{K}^n$. We denote by
$\mathbf{Val}$ the vector space of all continuous translation
invariant complex valued valuations and we write
$\mathbf{Val}_i^\pm$ for its subspace of all valuations of degree
$i$ and even/odd parity. An important result by McMullen
\textbf{\cite{McMullen77}} is that
\begin{equation} \label{eq_mcmullen}
\mathbf{Val} = \bigoplus_{0 \leq i \leq n} (\mathbf{Val}_i^+
\oplus \mathbf{Val}_i^-).
\end{equation}

In order to state our main theorem, we need the following basic
fact from the representation theory of the group $\mathrm{SO}(n)$:
The isomorphism classes of irreducible representations of
$\mathrm{SO}(n)$ are parametrized by their highest weights,
namely sequences of integers
$(\lambda_1,\lambda_2,\ldots,\lambda_{\lfloor n/2 \rfloor})$ such
that
\[\left \{\begin{array}{ll} \lambda_1 \geq \lambda_2 \geq \ldots \geq \lambda_{\lfloor n/2 \rfloor} \geq 0 &
\quad \mbox{for odd }n, \\
\lambda_1 \geq \lambda_2 \geq \ldots \geq \lambda_{n/2-1} \geq |\lambda_{n/2}| & \quad
\mbox{for even }n.
\end{array} \right .\]
(See Section 3 for the background material from representation
theory.)

The natural action of the group $\mathrm{SO}(n)$ on the space
$\mathbf{Val}$ is given by
\begin{displaymath}
 (\vartheta \phi)(K) = \phi(\vartheta ^{-1}K), \qquad \vartheta \in \mathrm{SO}(n),\, \phi \in \mathbf{Val}.
\end{displaymath}

Our main theorem is the following decomposition of the space
$\mathbf{Val}$ into irreducible $\mathrm{SO}(n)$-modules.

\begin{satz} \label{thm_decomposition}
Let $0 \leq i \leq n$. The space $\mathbf{Val}_i$ is the direct
sum of the \linebreak irreducible representations of
$\mathrm{SO}(n)$ with highest weights
$(\lambda_1,\ldots,\lambda_{\lfloor n/2 \rfloor})$ \linebreak
precisely satisfying the following additional conditions:
\begin{enumerate}
\item[(i)] $\lambda_j = 0$ for $j > \min\{i,n-i\}$;
\vspace{-0.1cm}
\item[(ii)] $|\lambda_j| \neq 1$ for $1 \leq j \leq \lfloor n/2 \rfloor$;
\vspace{-0.1cm}
\item[(iii)] $|\lambda_2| \leq 2$.
\end{enumerate}
In particular, under the action of $\mathrm{SO}(n)$ the space
$\mathbf{Val}_i$ is multiplicity free.
\end{satz}

\pagebreak

Earlier versions of Theorem \ref{thm_decomposition} for {\it even}
valuations were obtained in \textbf{\cite{Alesker01}} and
\textbf{\cite{AlBern}}. These results were subsequently applied
in the construction of new algebraic structures on the space
$\mathbf{Val}$ (see \textbf{\cite{Alesker04b, bernigfu06}}) which
provided the means for a fuller understanding of the integral
geometry of compact groups acting transitively on the unit sphere
(see e.g. \textbf{\cite{Alesker03, Bernig09, bernigfu10, fu06}}).

For the proof of Theorem \ref{thm_decomposition} we draw on
methods from representation \linebreak theory, differential
geometry and geometric measure theory. To be more specific, we use
a representation of smooth translation invariant valuations via
integral currents first obtained in \textbf{\cite{Alesker06a}}
and later refined in \textbf{\cite{Bernig07b}} and
\textbf{\cite{bernig10}} as well as an analysis of the action of
$\mathrm{SO}(n)$ on the space of translation invariant
differential forms on a contact manifold (see Sections 4 and 5).

\vspace{0.1cm}

We now state a reformulation of Theorem \ref{thm_decomposition} in
the language of tensor \linebreak valuations. Let
$(\Gamma,\varrho$) be a (finite dimensional, complex)
representation of $\mathrm{SO}(n)$. A continuous translation
invariant valuation with values in $\Gamma$ is called \linebreak
{\it $\mathrm{SO}(n)$ equivariant} if
\begin{displaymath}
 \phi(\vartheta K)=\varrho(\vartheta) \phi(K)
\end{displaymath}
for all $\vartheta \in \mathrm{SO}(n)$ and $K \in \mathcal{K}^n$.

\vspace{0.4cm}

\noindent {\bf Theorem $\mathbf{1'}$.\ } \emph{Let
$(\Gamma,\varrho)$ be an irreducible $\mathrm{SO}(n)$
representation and let $0 \leq i \leq n$. There exists a
non-trivial continuous translation invariant and $\mathrm{SO}(n)$
equivariant valuation of degree $i$ with values in $\Gamma$ if
and only if the highest weight of $\Gamma$ satisfies the
conditions (i)-(iii) from Theorem \ref{thm_decomposition}.
\linebreak This valuation is unique up to scaling.}

\vspace{0.4cm}

Since a finite dimensional representation of $\mathrm{SO}(n)$ can
be decomposed into a sum of irreducible representations, Theorem
$1'$ can be used to study the space of equivariant $\Gamma$-valued
valuations also for reducible $\Gamma$ (compare the examples in
Section 5).

The case of {\it symmetric} tensors, namely
$\Gamma=\mathrm{Sym}^kV$, has been intensively treated in
\textbf{\cite{hug_schneider_schuster_a, hug_schneider_schuster_b,
Ludwig:matrix, Alesker99, McMullen97}}. In these papers,
translation invariance is replaced by the more general {\it
isometry covariance}. In the recent article
\textbf{\cite{hug_schneider_schuster_b}}, Hug, Schneider and
R.~Schuster explicitly determined the dimension of the space of
all continuous isometry covariant tensor valuations of a fixed
rank and of a given degree of homogeneity. However, these
computations do not seem to give a basis of the subspace of
translation invariant tensor valuations. For the general,
non-symmetric, case, not much seems to be known except the
construction of $\Lambda^kV \otimes \Lambda^kV$-valued
translation invariant valuations in \textbf{\cite{Bernig06}}.

\pagebreak

\noindent {\bf Definition.} A map $\varphi: \mathcal{K}^n \times
\mathcal{K}^n \rightarrow \mathbb{C}$ is called a
\emph{bivaluation} if $\varphi$ is additive in each argument. A
bivaluation $\varphi$ is called {\it translation biinvariant} if
$\varphi$ is invariant under independent translations of its
arguments and $\varphi$ is said to have {\it bidegree} $(i,j)$ if
$\varphi(t K,s L)=t^is^j\varphi(K,L)$ for all $K,L \in
\mathcal{K}^n$ and $t, s > 0$. We say $\varphi$ is
\emph{$\mathrm{O}(n)$ invariant} (resp.\ $\mathrm{SO}(n)$
invariant) if $\varphi(\vartheta K,\vartheta L) = \varphi(K,L)$
for all $K, L \in \mathcal{K}^n$ and $\vartheta \in \mathrm{O}(n)$
(resp.\ $\vartheta \in \mathrm{SO}(n)$).

\vspace{0.3cm}

In their book on geometric probability, Klain and Rota
\textbf{\cite{Klain:Rota}} pose the \linebreak problem to classify
all ''invariant'' bivaluations. First such classification
\linebreak results were obtained recently by Ludwig
\textbf{\cite{Ludwig10a}}. In Section 6 we obtain a description
of all continuous translation biinvariant bivaluations which can
be seen as a starting point for systematic investigations of this
problem.

As an application of Theorem \ref{thm_decomposition}, we obtain
the following important \linebreak symmetry property of rigid
motion invariant homogeneous bivaluations.

\vspace{0.05cm}

\begin{satz} \label{thm_symmetry} If $\varphi: \mathcal{K}^n \times \mathcal{K}^n
\rightarrow \mathbb{R}$ is a continuous translation biinvariant
and $\mathrm{O}(n)$ invariant bivaluation of bidegree $(i,i)$, $0
\leq i \leq n$, then
\begin{equation} \label{eq_symmetry}
 \varphi(K,L)=\varphi(L,K)
\end{equation}
for every $K, L \in \mathcal{K}^n$.
\end{satz}

\vspace{0.05cm}

As a byproduct of our proof of Theorem \ref{thm_symmetry}, we also
obtain that if the bivaluation $\varphi$ is as above but merely
$\mathrm{SO}(n)$ invariant, then \eqref{eq_symmetry} still holds
true if $(i,n) \neq (2k+1,4k+2)$, $k \in \mathbb{N}$. If $n \equiv
2\! \mod 4$, then there exist \linebreak continuous translation
biinvariant and $\mathrm{SO}(n)$ invariant bivaluations of
\linebreak bidegree $\left (\frac{n}{2},\frac{n}{2} \right )$
which are \emph{not symmetric}.

\vspace{0.1cm}

The symmetry property established in Theorem 1 in combination with
techniques developed by Lutwak \textbf{\cite{lutwak85, lutwak93}}
can be used to obtain geometric inequalities for Minkowski
valuations. Recall that a map $\Phi: \mathcal{K}^n \rightarrow
\mathcal{K}^n$ is called a {\it Minkowski valuation} if $\Phi$ is
additive with respect to the usual Minkowski addition of convex
sets. We denote by $\mathbf{MVal}_i$ the set of all continuous
translation invariant Minkowski valuations of degree $i$.

A convex body $K$ is uniquely determined by its support function
$h(K,u)=\max\{u\cdot x: x \in K\}$, for $u \in S^{n-1}$. Among the
most important examples of Minkowski valuations is the projection
operator $\Pi \in \mathbf{MVal}_{n-1}$: The {\it projection body}
$\Pi K$ of $K$ is the convex body defined by
\[h(\Pi K,u)=\mathrm{vol}_{n-1}(K|u^{\bot}), \qquad u \in S^{n-1},  \]
where $K|u^{\bot}$ denotes the projection of $K$ onto the
hyperplane orthogonal to $u$. \linebreak For the special role of
the map $\Pi$ in the theory of valuations we refer to
\textbf{\cite{ludwig02}}.

\pagebreak

Let $0 \leq i \leq n - 1$. If $\Phi_i \in \mathbf{MVal}_i$ is
$\mathrm{O}(n)$ {\it equivariant}, i.e.\ $\Phi_i(\vartheta K)=
\vartheta \Phi_i K$ for all $K \in \mathcal{K}^n$ and $\vartheta
\in \mathrm{O}(n)$, then the map
\[\varphi(K,L)=V(\Phi_i K,L[i],B[n-i-1]), \qquad K, L \in \mathcal{K}^n,   \]
where $V(\Phi_i K,L[i],B[n-i-1])$ denotes the mixed volume of
$\Phi_i K$, $i$ copies of $L$ and $n - i - 1$ copies of the
Euclidean unit ball $B$, is a translation biinvariant and
$\mathrm{O}(n)$ invariant bivaluation of bidegree $(i,i)$. By
Theorem \ref{thm_symmetry}, it is symmetric in $K$ and $L$.

In the particular case where $i = n - 1$ and $\Phi_{n-1}=\Pi$,
this symmetry \linebreak property is well known. Its variants and
generalizations have been used \linebreak extensively for
establishing geometric inequalities related to convex and star
body valued valuations (see \textbf{\cite{goodeyzhang,
grinberg:zhang, habschu09, lutwak85, lutwak86, lutwak90,
lutwak93, LYZ2000a, LYZ2000b, LZ1997, Schu09}} and Section 7).
\linebreak Complex versions of the projection body were recently
studied in \textbf{\cite{Abardia_Bernig11}}, they satisfy similar
symmetry properties.

\vspace{0.1cm}

In the following we give one example of the type of inequalities
that can be derived from Theorem \ref{thm_symmetry}. To this end
let us recall a version of the classical Brunn--Minkowski
inequality. For $i \in \{0, \ldots, n\}$, let $V_i(K)$ denote the
$i$-th intrinsic volume of $K \in \mathcal{K}^n$. The
Brunn--Minkowski inequality for intrinsic volumes states the
following: If $2 \leq i \leq n$ and $K, L \in \mathcal{K}^n$ have
non-empty interior, then
\begin{equation} \label{intvolbm}
V_i(K+L)^{1/i} \geq V_i(K)^{1/i}+V_i(L)^{1/i},
\end{equation}
with equality if and only if $K$ and $L$ are homothetic.

In \textbf{\cite{lutwak85, lutwak93}} Lutwak obtained
inequalities of Brunn--Minkowski type for a well known family of
Minkowski valuations derived from the projection body operator.
As an application of Theorem \ref{thm_symmetry}, we show that
inequalities (\ref{intvolbm}) and Lutwak's inequalities for
derived projection operators of order $i$ are in fact part of a
larger family of Brunn--Minkowski type inequalities which hold
for all continuous translation invariant and $\mathrm{SO}(n)$
equivariant Minkowski valuations of a given degree.

\begin{satz} \label{thm_brunn_minkowski} Suppose that $\Phi_i \in \mathbf{MVal}_i$, $1 \leq i \leq n - 1$, is
$\mathrm{SO}(n)$ equivariant.
If $K, L \in \mathcal{K}^n$ have non-empty interior, then
\[V_{i+1}(\Phi_i(K + L))^{1/i(i+1)} \geq V_{i+1}(\Phi_iK)^{1/i(i+1)}+V_{i+1}(\Phi_iL)^{1/i(i+1)}.  \]
If $i \geq 2$ and $\Phi_i$ maps convex bodies with non-empty
interior to bodies with non-empty interior, then equality holds
if and only if $K$ and $L$ are homothetic.
\end{satz}

The special case of Theorem \ref{thm_brunn_minkowski} for {\it
even} Minkowski valuations was recently established by other
methods by the third author \textbf{\cite{Schu09}}.

\pagebreak

\centerline{\large{\bf{ \setcounter{abschnitt}{2}
\arabic{abschnitt}. Translation invariant valuations}}}

\vspace{0.6cm} \reseteqn \alpheqn \setcounter{theorem}{0}

In the following we collect background results on translation
invariant complex valued valuations needed in subsequent
sections. In particular, we recall the definition of
$\mathrm{O}(n)$ finite valuations and smooth valuations as well as
their representation via integral currents.

A classical theorem of Minkowski states that the volume of a
Minkowski linear combination $t_1K_1 + \ldots + t_k K_k$ of
convex bodies $K_1, \ldots, K_k$ can be expressed as a homogeneous
polynomial of degree $n$:
\begin{equation*}
V_n(t_1K_1 + \ldots +t_k K_k)=\sum \limits_{i_1,\ldots, i_n=1}^k
V(K_{i_1},\ldots,K_{i_n})t_{i_1}\cdots t_{i_n}.
\end{equation*}
The coefficients $V(K_{i_1},\ldots,K_{i_n})$ are called {\it mixed
volumes} of $K_{i_1}, \ldots, K_{i_n}$. Clearly,
$V(K,\ldots,K)=V_n(K)$. Moreover, mixed volumes are symmetric,
non-negative and multilinear with respect to Minkowski linear
combinations. They are also continuous with respect to the
Hausdorff metric and satisfy the following two properties:
\begin{itemize}
\item If $K, L \in \mathcal{K}^n$ such that $K \cup L \in \mathcal{K}^n$, and $\mathbf{C}=(K_1,\ldots,K_i)$, then
\[V_i(K,\mathbf{C})+V_i(L,\mathbf{C})=V_i(K \cup L,\mathbf{C})+V_i(K \cap L,\mathbf{C}),  \]
where $V_i(K,\mathbf{C})$ denotes the mixed volume
$V(K,\ldots,K,K_1,\ldots,K_i)$.
\item Mixed volumes are invariant under independent translations of their arguments
and they are invariant under simultaneous unimodular linear
transformations, i.e., if $K_1, \ldots, K_n \in \mathcal{K}^n$
and $A \in \mathrm{SL}(n)$, then
\[V(AK_1,\ldots,AK_n) = V(K_1,\ldots,K_n).   \]
\end{itemize}

Recall that we denote by $\mathbf{Val}$ the vector space of
continuous translation invariant {\it complex valued} valuations
and we write $\mathbf{Val}_i^\pm$ for its subspaces of all
valuations of degree $i$ and even/odd parity.

It is easy to see that the space $\mathbf{Val}_0$ is
one-dimensional. The analogous (non-trivial) statement for
$\mathbf{Val}_n$ was proved by Hadwiger
\textbf{\cite[\textnormal{p.\ 79}]{hadwiger51}}.

The following consequence of McMullen's decomposition \eqref{eq_mcmullen} is well
known.

\begin{koro} \label{normval} Let $C \in \mathcal{K}^n$ be a fixed convex body with
non-empty interior. The space $\mathbf{Val}$ becomes a Banach
space under the norm
\[\|\phi \| = \sup \{|\phi(K)|: K \subseteq C \}.  \]
Moreover, a different choice of $C$ yields an equivalent norm.
\end{koro}

The natural continuous action of the group $\mathrm{GL}(n)$ on
the Banach space $\mathbf{Val}$ is given by
\[(A \phi)(K)=\phi(A^{-1}K), \qquad A \in \mathrm{GL}(n),\, \phi \in \mathbf{Val}. \]
Note that the subspaces $\mathbf{Val}_i^{\pm} \subseteq
\mathbf{Val}$ are invariant under this $\mathrm{GL}(n)$ action.

The following result is known as the {\it Irreducibility
Theorem}. It implies a conjecture by McMullen that the linear
combinations of mixed volumes form a dense subspace in
$\mathbf{Val}$.

\begin{theorem} \label{irr} \emph{(\!\!\textbf{\cite{Alesker01}})} The
natural action of $\mathrm{GL}(n)$ on $\mathbf{Val}_i^{\pm}$ is
irreducible for every $i \in \{0, \ldots, n\}$.
\end{theorem}

In the following it will be important for us to work with two different dense
subsets of valuations in $\mathbf{Val}$:

\vspace{0.3cm}

\noindent {\bf Definition} \emph{A valuation $\phi \in
\mathbf{Val}$ is called $\mathrm{O}(n)$ finite if the
$\mathrm{O}(n)$ orbit of $\phi$, i.e.\ the subspace
$\mathrm{span}\{\vartheta \phi: \vartheta \in \mathrm{O}(n)\}$, is
finite dimensional.}

\emph{A valuation $\phi \in
\mathbf{Val}$ is called smooth if the map $\mathrm{GL}(n)
\rightarrow \mathbf{Val}$ defined by $A \mapsto A\phi$ is
infinitely differentiable.}

\vspace{0.3cm}

The notions of $\mathrm{O}(n)$ finite and smooth valuations are
special cases of more general well known concepts from
representation theory (see e.g.\ \textbf{\cite{wallach1}}).

We denote the space of continuous translation invariant and
$\mathrm{O}(n)$ finite valuations by $\mathbf{Val}^{f}$ and we
write $\mathbf{Val}^{\infty}$ for the space of smooth translation
invariant valuations. For the subspaces of homogeneous valuations
of given parity in $\mathbf{Val}^{f}$ and $\mathbf{Val}^{\infty}$
we write $\mathbf{Val}_i^{\pm,f}$ and
$\mathbf{Val}_i^{\pm,\infty}$, respectively.

It is well known (cf.\ \textbf{\cite[\textnormal{p.\
141}]{broecker_tomdieck}} and \textbf{\cite[\textnormal{p.\
32}]{wallach1}}) that the set of $\mathrm{O}(n)$ finite valuations
$\mathbf{Val}_i^{\pm,f}$ is a dense $\mathrm{O}(n)$ invariant
subspace of $\mathbf{Val}_i^{\pm}$ and that the set of smooth
valuations $\mathbf{Val}_i^{\pm,\infty}$ is a dense
$\mathrm{GL}(n)$ invariant subspace of $\mathbf{Val}_i^{\pm}$.
Moreover, $\mathbf{Val}^{f} \subseteq \mathbf{Val}^{\infty}$ and
from \eqref{eq_mcmullen} one easily deduces that the spaces
$\mathbf{Val}^{f}$ and $\mathbf{Val}^{\infty}$ admit direct sum
decompositions into their corresponding subspaces of homogeneous
valuations of given parity.

\vspace{0.2cm}

An equivalent description of smooth valuations can be given in
terms of the normal cycle map. Let $SV = V \times S^{n-1}$ denote
the unit sphere bundle on $V$. For $K \in \mathcal{K}^n$ and $x
\in \partial K$, we write $N(K,x)$ for the normal cone of $K$ at
$x$. The {\it normal cycle} (or generalized normal bundle) of a
convex body $K$ is the Lipschitz submanifold of $SV$ defined by
\[\mathrm{nc}(K) = \{(x,u) \in SV: x \in \partial K, u \in N(K,x)\}. \]

For $0 \leq i \leq n - 1$, let $\Omega^{i,n-i-1}$ denote the
space of smooth translation invariant differential forms of
bidegree $(i,n-i-1)$ on $SV$. The following result is a special
case of \textbf{\cite[\textnormal{Theorem 5.2.1}]{Alesker06a}}:

\begin{lem} \label{smviacurr}
If $0 \leq i \leq n - 1$, then the map $\nu: \Omega^{i,n-i-1}
\rightarrow \mathbf{Val}_i^{\infty}$, defined \nolinebreak by
\begin{equation} \label{eq_normal_cycle}
\nu(\omega)(K)=\int_{\mathrm{nc}(K)} \omega,
\end{equation}
is surjective.
\end{lem}

The kernel of the map $\nu$ was described in
\textbf{\cite{Bernig07b}} in terms of the Rumin operator
\textbf{\cite{rumin}}, a second order differential operator which
acts on smooth forms on the sphere bundle. A refined version of
this result (stated in Section 4 as Theorem \ref{prop_exact_seq})
was recently proved in \textbf{\cite{bernig10}} and will be
crucial in the proof of Theorem \ref{thm_decomposition}. We also remark that recently a
broader notion of smooth valuations in the setting of smooth
manifolds was introduced, see \textbf{\cite{Alesker06a}}. The
classical concept of valuations as used in this article is in
some sense an infinitesimal version of this more general notion.

The description of smooth valuations provided by Lemma
\ref{smviacurr} was the main tool used in
\textbf{\cite{Bernig07b}} to establish a {\it Hard Lefschetz
Theorem} for translation invariant valuations (see also
\textbf{\cite{Alesker03, alesker10}}). The next statement is an
immediate consequence of this result:

\begin{theorem} \label{hardlef}
For every $i \in \{0, \ldots, n\}$, the spaces
$\mathbf{Val}_i^{\infty}$ and $\mathbf{Val}_{n-i}^{\infty}$ are
isomorphic as $\mathrm{SO}(n)$ modules.
\end{theorem}

\vspace{1cm}

\centerline{\large{\bf{ \setcounter{abschnitt}{3}
\arabic{abschnitt}. Irreducible representations of
$\mathrm{SO}(n)$ and $\mathrm{O}(n)$}}}

\vspace{0.6cm} \reseteqn \alpheqn \setcounter{theorem}{0}

In this section we recall some well known results concerning
irreducible representations of the groups $\mathrm{SO}(n)$ and
$\mathrm{O}(n)$, $n \geq 3$. As a general reference for this
material we recommend the books by Br\"ocker and tom Dieck
\textbf{\cite{broecker_tomdieck}}, Fulton and Harris
\textbf{\cite{fultonharris}}, and Goodman and Wallach
\textbf{\cite{goodwallach}}.

\pagebreak

Since $\mathrm{SO}(n)$ and $\mathrm{O}(n)$ are compact Lie groups,
all their irreducible \linebreak representations are finite
dimensional. The equivalence classes of irreducible {\it complex}
representations of $\mathrm{SO}(n)$ are indexed by their highest
weights, namely $\lfloor n/2 \rfloor$-tuples of integers
$(\lambda_1,\lambda_2,\ldots,\lambda_{\lfloor n/2 \rfloor})$ such
that
\begin{equation} \label{heiwei}
\left \{\begin{array}{ll} \lambda_1 \geq \lambda_2 \geq \ldots \geq \lambda_{\lfloor n/2 \rfloor} \geq 0 & \quad \mbox{for odd }n, \\
\lambda_1 \geq \lambda_2 \geq \ldots \geq \lambda_{n/2-1} \geq
|\lambda_{n/2}| & \quad \mbox{for even }n.
\end{array} \right .
\end{equation}

We refer to \cite{broecker_tomdieck} or \cite{adams} for an introduction to highest weights. 

In the following we use $\Gamma_\lambda$ to denote any isomorphic
copy of an irreducible representation of $\mathrm{SO}(n)$ with
highest weight
$\lambda=(\lambda_1,\lambda_2,\ldots,\lambda_{\lfloor n/2
\rfloor})$.

\vspace{0.4cm}

\noindent {\bf Examples:}

\begin{enumerate}
\item[(a)] The only one dimensional (complex) representation of
$\mathrm{SO}(n)$ is the trivial representation; it corresponds to
the $\mathrm{SO}(n)$ module $\Gamma_{(0,\ldots,0)}$.

\item[(b)] We denote by $V_{\mathbb{C}} = V \otimes \mathbb{C}$ the complexification of
$V$. The standard representation of $\mathrm{SO}(n)$ on
$V_{\mathbb{C}}$ corresponds to $\Gamma_{(1,0,\ldots,0)}$.

\item[(c)] For every $0 \leq i \leq \lfloor n/2 \rfloor -1$, the exterior power $\Lambda^i V_{\mathbb{C}}$ is
an irreducible $\mathrm{SO}(n)$ module with
$\lambda=(1,\ldots,1,0,\ldots,0)$, where $1$ appears $i$ times.

If $n = 2k + 1$ is odd, the exterior power $\Lambda^k
V_{\mathbb{C}}$ is also irreducible; but if $n = 2k$ is even, it
splits as $\Lambda^k V_{\mathbb{C}} = \Gamma_{(1,\ldots,1)}\oplus
\Gamma_{(1,\ldots,1,-1)}$

For every $i \in \{0,\ldots,n\}$, there is a natural isomorphism
of $\mathrm{SO}(n)$ \linebreak modules
\begin{equation} \label{lamknmink}
 \Lambda^i V_{\mathbb{C}} \cong \Lambda^{n-i} V_{\mathbb{C}}.
\end{equation}
\item[(d)] For $k \geq 2$, the symmetric power $\mathrm{Sym}^k V_{\mathbb{C}}$
is not irreducible as $\mathrm{SO}(n)$ module; its decomposition
into irreducible submodules is given by
\begin{equation} \label{eq_decomp_symmetric}
 \mathrm{Sym}^k V_{\mathbb{C}} =\bigoplus_{j=0}^{\lfloor k/2 \rfloor}
 \Gamma_{(k-2j,0,\ldots,0)}.
\end{equation}
\end{enumerate}

A description of the irreducible representations of the full
orthogonal group $\mathrm{O}(n)$ can be given in terms of the
irreducible representations of its identity component
$\mathrm{SO}(n)$ (cf.\ \textbf{\cite[\textnormal{p.\
249}]{goodwallach}}). The main difference arises from the fact
that $\mathrm{O}(n)$ has a non-trivial one dimensional
representation, called the {\it determinant representation},
which corresponds to the $\mathrm{O}(n)$ module
$\Lambda^nV_{\mathbb{C}}$.

\begin{lem} \label{onirr} Let $\lambda =
(\lambda_1,\ldots,\lambda_{\lfloor n/2 \rfloor})$ be a tuple of
integers satisfying (\ref{heiwei}).
\begin{enumerate}
\item[(a)] If $n$ is odd, then the irreducible representation
$\Gamma_\lambda$ of $\mathrm{SO}(n)$ is the \linebreak restriction
of two non-isomorphic irreducible $\mathrm{O}(n)$ representations
$\bar{\Gamma}_\lambda$ and $\bar{\Gamma}_\lambda \otimes
\Lambda^nV_{\mathbb{C}}$.
\item[(b)] If $n$ is even and $\lambda_{n/2} = 0$, then the irreducible representation
$\Gamma_\lambda$ \linebreak of $\mathrm{SO}(n)$ is the restriction
of two non-isomorphic irreducible $\mathrm{O}(n)$ \linebreak
representations $\bar{\Gamma}_\lambda$ and $\bar{\Gamma}_\lambda
\otimes \Lambda^nV_{\mathbb{C}}$. If $\lambda_{n/2} \neq 0$, then
$\Gamma_\lambda$ is not such a restriction.
\item[(c)] If $n$ is even and $\lambda_{n/2} \neq 0$, then the
$\mathrm{SO}(n)$ representation $\Gamma_\lambda \oplus
\Gamma_{\lambda'}$, where
$\lambda':=(\lambda_1,\ldots,\lambda_{n/2-1},-\lambda_{n/2})$, is
the restriction of an irreducible  $\mathrm{O}(n)$ representation
$\bar{\Gamma}_\lambda$.
\end{enumerate}
Moreover, all irreducible representations of $\mathrm{O}(n)$ are
determined in this way.
\end{lem}

Let $\Gamma$ be a (not necessarily irreducible) complex
$\mathrm{SO}(n)$ or $\mathrm{O}(n)$ module. Recall that the dual
representation is defined on the dual space $\Gamma^*$ by
\[(\vartheta\,u^*)(v) = u^*(\vartheta^{-1} v), \qquad \vartheta \in \mathrm{SO}(n), u^* \in \Gamma^*, v \in \Gamma.  \]
We say that $\Gamma$ is {\it self-dual} if $\Gamma$ and $\Gamma^*$
are isomorphic representations. The module $\Gamma$ is called {\it
real} if there exists a non-degenerate symmetric $\mathrm{SO}(n)$
invariant, or $\mathrm{O}(n)$ respectively, bilinear form on
$\Gamma$. In particular, if $\Gamma$ is
real, then $\Gamma$ is also self-dual.

The following lemma (cf.\ \textbf{\cite[\textnormal{p.\
292}]{broecker_tomdieck}})  will be critical in the proof of
Theorem \ref{thm_symmetry}:

\begin{lem} \label{lemma_real_reps} Let $\lambda =
(\lambda_1,\ldots,\lambda_{\lfloor n/2 \rfloor})$ be a tuple of
integers satisfying (\ref{heiwei}).
\begin{enumerate}
\item[(a)] If $n \not \equiv 2 \mod 4$, then all representations of $\mathrm{SO}(n)$ are
real.
\item[(b)] If $n \equiv 2 \mod 4$, then the irreducible representation $\Gamma_\lambda$ of $\mathrm{SO}(n)$ is real if
and only if $\lambda_{n/2} = 0$. If $\lambda_{n/2} \neq 0$, then
the dual of $\Gamma_\lambda$ is $\Gamma_{\lambda'}$.
\end{enumerate}
Moreover, all representations of $\mathrm{O}(n)$ are real.
\end{lem}

An essential tool in the classification of irreducible modules of
a compact group is the character of a representation: Let $\Gamma$
be a finite dimensional (complex) $\mathrm{SO}(n)$ module and let
$\varrho: \mathrm{SO}(n) \rightarrow \mathrm{GL}(\Gamma)$ be the
corresponding representation. The {\it character} of $\Gamma$ is
the function  $\mathrm{char}\, \Gamma: \mathrm{SO}(n) \rightarrow
\mathbb{C}$ defined by
\[(\mathrm{char}\, \Gamma)(\vartheta) = \mathrm{Tr}\,\varrho(\vartheta), \]
where $\mathrm{Tr}\,\varrho(\vartheta)$ is the trace of the linear
map $\varrho(\vartheta): \Gamma \rightarrow \Gamma$.

A complex representation is determined up to isomorphism by its
\linebreak character. Moreover, from properties of the trace map,
one immediately obtains several useful arithmetic properties of
characters: If $\Gamma$ and $\Theta$ are finite dimensional
$\mathrm{SO}(n)$ modules, then
\begin{equation} \label{propchar1}
\mathrm{char}(\Gamma \oplus
\Theta)=\mathrm{char}\,\Gamma+\mathrm{char}\,\Theta
\end{equation}
and
\begin{equation} \label{propchar2}
\mathrm{char}(\Gamma \otimes \Theta)=\mathrm{char}\,\Gamma \cdot
\mathrm{char}\,\Theta.
\end{equation}

The character of the irreducible $\mathrm{SO}(n)$ modules
$\Gamma_\lambda$ with highest weights $\lambda =
(\lambda_1,\ldots,\lambda_{\lfloor n/2 \rfloor})$ are described
by Weyl's character formula. However, more important for us is a
consequence of this description, known as the {\it second
determinantal formula}, which we describe in the following.

Let $\lambda = (\lambda_1,\ldots,\lambda_{\lfloor n/2 \rfloor})$
be a tuple of {\it non-negative} integers satisfying
(\ref{heiwei}). \linebreak We define the $\mathrm{SO}(n)$ module
$\bar{\Gamma}_\lambda$ by
\begin{equation*}
\bar{\Gamma}_\lambda := \left \{\begin{array}{ll}
\Gamma_{\lambda} \oplus \Gamma_{\lambda'} & \quad
\mbox{if $n$ is even and } \lambda_{n/2} \neq 0, \\
\Gamma_{\lambda}  & \quad \mbox{otherwise.}
\end{array} \right .
\end{equation*}
The second determinantal formula expresses
$\mathrm{char}\,\bar{\Gamma}_\lambda$ as a polynomial in the
characters $E_i$ of the fundamental representations
$\Lambda^iV_{\mathbb{C}}$, $i \in \mathbb{Z}$. (Note that
$E_0=E_n=1$ and that we use the convention $E_i=0$ for $i<0$ or
$i>n$.)

Given a tuple of {\it non-negative} integers $\lambda =
(\lambda_1,\ldots,\lambda_{\lfloor n/2 \rfloor})$ satisfying
(\ref{heiwei}), recall that the {\it conjugate} of $\lambda$ is
the $s := \lambda_1$ tuple $\mu = (\mu_1,\ldots,\mu_s)$ defined by
saying that $\mu_j$ is the number of terms in $\lambda$ that are
greater than or equal $j$. The second determinantal formula (cf.\
\textbf{\cite[\textnormal{p.\ 409}]{fultonharris}}) can be stated
as follows:

\begin{theorem} \label{secdet} Let $\lambda =
(\lambda_1,\ldots,\lambda_{\lfloor n/2 \rfloor})$ be a tuple of
non-negative integers satisfying (\ref{heiwei}) and let
$\mu=(\mu_1,\ldots,\mu_s)$ be the conjugate of $\lambda$. The
character of $\bar{\Gamma}_\lambda$ equals the determinant of the
$s \times s$-matrix whose $i$-th row is given by
\begin{equation} \label{secdetmat}
\left(
\begin{array}{c c c c} E_{\mu_i-i+1} & E_{\mu_i-i+2}+E_{\mu_i-i}
& \cdots & E_{\mu_i-i+s}+E_{\mu_i-i-s+2}\end{array}\right).
\end{equation}
\end{theorem}

It is sometimes convenient for us to take $s > \lambda_1$ in the
definition of the conjugate of $\lambda$. This just introduces
additional zeros at the end of the conjugate tuple. However, note
that this does not change the determinant of the matrix defined
by (\ref{secdetmat}).

In the following we use $\#(\lambda,j)$ to denote the number of
terms in a \linebreak tuple of (non-negative) integers $\lambda =
(\lambda_1,\ldots,\lambda_{\lfloor n/2 \rfloor})$ which are equal
to $j$. As a consequence of Theorem \ref{secdet}, we note the
following auxiliary result which will be needed in the proof of
Theorem \ref{thm_decomposition}.

\begin{koro} \label{cor_clebsch_gordan}
If $i, j \in \mathbb{N}$ are such that $n/2 \leq i \leq n$ and $i
+ j \leq n$, then
\begin{equation} \label{clebgoreq1}
E_iE_j-E_{i-1}E_{j-1}=\sum_\lambda \mathrm{char}\,
\bar{\Gamma}_\lambda,
\end{equation}
where the sum ranges over all $\lfloor n/2 \rfloor$-tuples of
non-negative integers \linebreak $\lambda=
(\lambda_1,\ldots,\lambda_{\lfloor n/2 \rfloor})$ satisfying
(\ref{heiwei}) and
\begin{equation} \label{clebgoreq2}
\lambda_1 \leq 2, \qquad \#(\lambda,1)=n-i-j, \qquad
\#(\lambda,2) \leq j.
\end{equation}
\end{koro}

\noindent {\it Proof}\,: If $\lambda =
(\lambda_1,\ldots,\lambda_{\lfloor n/2 \rfloor})$ is a tuple of
non-negative integers satisfying (\ref{heiwei}) and
(\ref{clebgoreq2}), then the conjugate of $\lambda$ is given by
$\mu=(\mu_1,\mu_2)$, where $\mu_2 = \#(\lambda,2) \leq j$ and
$\mu_1-\mu_2 = \#(\lambda,1)=n-i-j$. Thus, by Theorem
\ref{secdet}, the character of $\bar{\Gamma}_{\lambda}$ is given
by
\begin{displaymath}
\mathrm{char}\,\bar{\Gamma}_\lambda = \det \left(
\begin{array}{c c} E_{\mu_2+k} & E_{\mu_2+k+1}+E_{\mu_2+k-1} \\
E_{\mu_2-1} & E_{\mu_2}+E_{\mu_2-2} \end{array}\right),
\end{displaymath}
where $k= n - i - j$. Consequently, the right hand side of
(\ref{clebgoreq1}) is
\begin{eqnarray*}
 \sum_\lambda \mathrm{char}\,\bar{\Gamma}_\lambda & =
 & \sum_{\mu_2=0}^j  \big ( E_{\mu_2+k}(E_{\mu_2}+E_{\mu_2-2})- E_{\mu_2-1}(E_{\mu_2+k+1}+E_{\mu_2+k-1})  \big )\\
 & = & E_{n-i}E_j-E_{n-(i-1)}E_{j-1}.
\end{eqnarray*}
To finish the proof, note that $E_{n-i}=E_i$ by
(\ref{lamknmink}). \hfill $\blacksquare$

\vspace{0.4cm}

An important class of (infinite dimensional) representations of a
Lie group $G$ are those induced from closed subgroups $H$ of $G$.
Although in this article we will only need the case $G =
\mathrm{SO}(n)$ and $H=\mathrm{SO}(n-1)$, we shall explain this
construction for a general compact Lie group $G$ and its closed
subgroup $H$. To this end, for any finite dimensional complex
vector space $\Gamma$, we denote by $C^{\infty}(G;\Gamma)$ the
space of all smooth functions $f:G\rightarrow \Gamma$.

\pagebreak

If $\Theta$ is any representation of $G$, clearly we obtain a
representation $\mathrm{Res}_H^G\Theta$ of $H$ by restriction.
Conversely, each $H$ module $\Gamma$ induces a representation of
$G$ as follows: Let $\mathrm{Ind}_H^G\Gamma \subseteq
C^{\infty}(G;\Gamma)$ be the space of functions defined by
\[\mathrm{Ind}_H^G \Gamma:=\left\{f \in C^{\infty}(G;\Gamma): f(gh)=h^{-1}f(g) \mbox{ for all } g \in G, h \in
H \right\}. \] The (smooth) {\it induced representation} of $G$ on
$\mathrm{Ind}_H^G \Gamma$ is now given by left translation
\[(gf)(u)=f(g^{-1}u), \qquad g,u \in G.\]

A basic result on induced representations is the well known
{\it Frobenius Reciprocity Theorem} (cf.\
\textbf{\cite[\textnormal{p.\ 523}]{goodwallach}}):

\begin{theorem} \label{frobrec}
If $\Theta$ is a $G$ module and $\Gamma$ is an $H$ module, then
there is a canonical vector space isomorphism
\[\mathrm{Hom}_G(\Theta,\mathrm{Ind}_H^G\Gamma) \cong \mathrm{Hom}_H(\mathrm{Res}_H^G
\Theta,\Gamma).\]
\end{theorem}

Here, $\mathrm{Hom}_G$ denotes the space of continuous linear $G$
equivariant maps.

Recall that if $\Theta$ is an irreducible $G$ module, by Schur's
lemma, the \linebreak {\it multiplicity} $m(\Xi,\Theta)$ of
$\Theta$ in an arbitrary $G$ module $\Xi$ is given by
\[m(\Xi,\Theta)=\dim \mathrm{Hom}_G(\Xi,\Theta)=\dim \mathrm{Hom}_G(\Theta,\Xi).  \]
Thus, by the Frobenius Reciprocity Theorem, if $\Theta$ and
$\Gamma$ are irreducible, then the multiplicity of $\Theta$ in
$\mathrm{Ind}_H^G\Gamma$ equals the multiplicity of $\Gamma$ in
$\mathrm{Res}_H^G \Theta$.

In order to apply Theorem \ref{frobrec} in our situation, where
$G=\mathrm{SO}(n)$ and $H=\mathrm{SO}(n-1)$, we will need a
formula for decomposing
$\mathrm{Res}_{\mathrm{SO}(n-1)}^{\mathrm{SO}(n)} \Gamma$ into
irreducible $\mathrm{SO}(n-1)$ modules. This is the content of
the following {\it branching theorem} (cf.\
\textbf{\cite[\textnormal{p.\ 426}]{fultonharris}}):

\begin{theorem} \label{thm_branching}
If $\Gamma_{\lambda}$, with $\lambda =
(\lambda_1,\ldots,\lambda_{\lfloor n/2 \rfloor})$ satisfying
(\ref{heiwei}), is an \linebreak irreducible representation of
$\mathrm{SO}(n)$, then
\begin{equation}
\mathrm{Res}_{\mathrm{SO}(n-1)}^{\mathrm{SO}(n)}
\Gamma_{\lambda}=\bigoplus_\mu \Gamma_\mu,
\end{equation}
where the sum ranges over all $\mu=(\mu_1,\ldots,\mu_k)$ with
$k:=\lfloor (n-1)/2 \rfloor$ and
\begin{equation*}
\left \{\begin{array}{ll} \lambda_1 \geq \mu_1\geq \lambda_2 \geq
\mu_2 \geq \ldots \geq \mu_{k-1} \geq \lambda_{\lfloor n/2 \rfloor} \geq |\mu_k| & \quad \mbox{for odd }n, \\
\lambda_1 \geq \mu_1 \geq \lambda_2 \geq \mu_2 \geq \ldots \geq
\mu_{k} \geq |\lambda_{n/2}| & \quad \mbox{for even }n.
\end{array} \right .
\end{equation*}
\end{theorem}

\pagebreak

\centerline{\large{\bf{ \setcounter{abschnitt}{4}
\arabic{abschnitt}. The Rumin--de Rham complex}}}

\vspace{0.7cm} \reseteqn \alpheqn \setcounter{theorem}{0}

We state in this section a refinement of the description of
translation invariant smooth valuations via integral currents. We
also establish an \linebreak auxiliary  result which will enable
us to subsequently employ the machinery from representation theory
explained in Section 3.

Recall that $SV = V \times S^{n-1}$ denotes the unit sphere
bundle. The natural smooth (left) action of $\mathrm{SO}(n)$ on
$SV$ is given by
\begin{equation} \label{action17}
l_{\vartheta}(x,u):=(\vartheta x,\vartheta u), \qquad \vartheta
\in \mathrm{SO}(n), (x,u) \in SV.
\end{equation}
Similarly, each $y \in V$ determines a smooth map $t_y: SV
\rightarrow SV$ by
\begin{equation} \label{trans17}
t_y(x,u)=(x+y,u), \qquad (x,u) \in SV.
\end{equation}

The canonical {\it contact form} $\alpha$ on $SV$ is the one form
defined by
\[\alpha|_{(x,u)}(w)= \langle u, d_{(x,u)}\pi(w)\rangle, \qquad w \in T_{(x,u)}SV,\]
where $\pi:SV \to V$ denotes the canonical projection and
$d_{(x,u)}\pi$ its differential at $(x,u) \in SV$. In this way,
$SV$ becomes a $2n - 1$ dimensional contact manifold. The kernel
of $\alpha$ defines the {\it contact distribution} $Q:=\ker
\alpha$. The restriction of $d\alpha$ to $Q$ is a
non-degenerate two form. In this way, each $Q_{(x,u)}$ becomes a symplectic vector space.

The {\it Reeb vector field} $R$ on $SV$ is defined by
$R_{(x,u)}=(u,0)$. It is the unique vector field on $SV$ such that
$\alpha(R)=1$ and $i_Rd\alpha = 0$, where $i_Rd\alpha$ denotes
the interior product of $R$ and $d\alpha$. At each point $(x,u)$,
$Q_{(x,u)}$ is the orthogonal sum of two copies of $T_u S^{n-1}$
and, consequently, we have an orthogonal splitting of the tangent
space $T_{(x,u)}SV$ given by
\begin{equation} \label{tangsplit}
T_{(x,u)}SV=\mathrm{span}_{\mathbb{R}}R_{(x,u)} \oplus T_u
S^{n-1} \oplus T_u S^{n-1}.
\end{equation}

The product structure of $SV$ induces a bigrading on the vector
space $\Omega^*(SV)$ of complex valued smooth differential forms
given by
\[\Omega^*(SV)=\bigoplus_{i,j} \Omega^{i,j}(SV),\]
where $\Omega^{i,j}(SV)$ denotes the subspace of $\Omega^*(SV)$ of
forms of bidegree $(i,j)$. We write $\Omega^{i,j} \subseteq
\Omega^{i,j}(SV)$ for the subspace of translation invariant
forms, i.e.,
\[\Omega^{i,j}=\{\omega \in \Omega^{i,j}(SV): t_y^*\omega = \omega \mbox{ for all } y \in V\}. \]
Here, $t_y^*$ is the pullback of the map $t_y: SV \rightarrow SV$
defined in (\ref{trans17}). Note that the restriction of the
exterior derivative $d$ to $\Omega^{i,j}$ has bidegree $(0,1)$.

The vector space $\Omega^{i,j}$ becomes an $\mathrm{SO}(n)$
module under the (continuous) action
\[\vartheta \omega = l_{\vartheta^{-1}}^* \omega, \qquad \vartheta \in \mathrm{SO}(n), \omega \in \Omega^{i,j}.   \]

An important $\mathrm{SO}(n)$ submodule of $\Omega^{i,j}$ is
given by the space $\Omega_v^{i,j}$ of {\it vertical} forms,
defined by
\[\Omega_v^{i,j} := \{\omega \in \Omega^{i,j}: \alpha \wedge \omega=0\}.  \]
Note that a differential form $\omega \in \Omega^{i,j}$ is
vertical if and only if it vanishes on the contact distribution
$Q$ of $SV$.

The $\mathrm{SO}(n)$ submodule $\Omega_h^{i,j} \subseteq
\Omega^{i,j}$ of {\it horizontal} forms, is given by
\begin{align*}
\Omega_h^{i,j} := \{\omega \in \Omega^{i,j}: i_R \omega=0\} \cong
\Omega^{i,j} / \Omega_v^{i,j}.
\end{align*}
It follows from (\ref{tangsplit}) and the definition of
$\Omega_h^{i,j}$ that $\omega \in \Omega^{i,j}$ is horizontal if
and only if
\[\omega|_{(x,u)} \in  \Lambda^i T_u^*S^{n-1}
\otimes \Lambda^j T_u^*S^{n-1} \otimes \mathbb{C}  \] for every $x
\in V$ and each $u \in S^{n-1}$. In the following we will
therefore simply write $\omega|_u$ instead of $\omega|_{(x,u)}$
whenever $\omega \in \Omega_h^{i,j}$ and $(x,u) \in SV$.

We now fix a point $u_0 \in S^{n-1}$ and let $\mathrm{SO}(n-1)$
be the stabilizer of $\mathrm{SO}(n)$ at $u_0$. For $u \in
S^{n-1}$, we denote by $W_u := T_uS^{n-1} \otimes \mathbb{C}$ the
complexification of the tangent space $T_uS^{n-1}$ and we write $W_0$ to denote $W_{u_0}$.

\begin{lem} \label{lemma_induced_reps}
For $i,j \in \mathbb{N}$, there is an isomorphism of
$\mathrm{SO}(n)$ modules
\begin{displaymath}
\Omega_h^{i,j} \cong
\mathrm{Ind}_{\mathrm{SO}(n-1)}^{\mathrm{SO}(n)} (\Lambda^i
W_0^* \otimes \Lambda^j W_0^*).
\end{displaymath}
\end{lem}

\noindent {\it Proof}\,: First note that, for each $\vartheta \in
\mathrm{SO}(n)$, the differential of the map $l_{\vartheta}: SV
\rightarrow SV$ defined in (\ref{action17}) induces a linear
isomorphism
\[\widehat{d_{u_0}l_{\vartheta}}:=(d_{u_0}l_\vartheta)^*: \Lambda^i
W_{\vartheta u_0}^* \otimes \Lambda^j W_{\vartheta u_0}^*
\rightarrow \Lambda^i W_{0}^* \otimes \Lambda^j W_{0}^*. \]
Moreover, the natural representation of the group
$\mathrm{SO}(n-1)$ on the space $\Lambda^i W_{0}^* \otimes
\Lambda^j W_{0}^*$ is given by $\eta \mapsto
\widehat{d_{u_0}l_{\eta^{-1}}}$.

Suppose now that $\omega \in \Omega_h^{i,j}$. We define
$f_{\omega}: \mathrm{SO}(n) \rightarrow \Lambda^i W_{0}^*
\otimes \Lambda^j W_{0}^*$ by
\[f_{\omega}(\vartheta) =\widehat{d_{u_0}l_{\vartheta}}(\omega|_{\vartheta
u_0}).\]
Clearly, we have $f_{\omega}(\vartheta \eta) = \eta^{-1}
f_{\omega}(\vartheta)$ for every $\vartheta \in \mathrm{SO}(n)$
and $\eta \in \mathrm{SO}(n-1)$. This shows that $f_{\omega} \in
\mathrm{Ind}_{\mathrm{SO}(n-1)}^{\mathrm{SO}(n)} (\Lambda^i
W_{0}^* \otimes \Lambda^j W_{0}^*)$.

Conversely, let $f \in
\mathrm{Ind}_{\mathrm{SO}(n-1)}^{\mathrm{SO}(n)} (\Lambda^i
W_{0}^* \otimes \Lambda^j W_{0}^*)$. We define a horizontal
form $\omega_f \in \Omega^{i,j}_h$ by
\[\omega_f|_{\vartheta u_0} = \widehat{d_{u_0}l_{\vartheta}}^{-1}(
f(\vartheta)).\] It is not difficult to show that $\omega$ is
well defined, i.e.\ if $\vartheta u_0 = \vartheta' u_0$ for some
$\vartheta, \vartheta' \in \mathrm{SO}(n)$, then
$\omega|_{\vartheta u_0}=\omega|_{\vartheta' u_0}$.

The observation that the $\mathrm{SO}(n)$ equivariant linear maps $\omega \mapsto f_{\omega}$
and $f \mapsto \omega_f$ are inverse to each other finishes the
proof. \hfill $\blacksquare$

\vspace{0.4cm}

Let $\mathcal{I}^{i,j}$ denote the $\mathrm{SO}(n)$ invariant
subspace of $\Omega^{i,j}$ defined by
\[\mathcal{I}^{i,j} :=\{ \omega \in \Omega^{i,j}: \omega= \alpha
\wedge \xi+ d\alpha \wedge \psi,\, \xi \in \Omega^{i-1,j}, \psi
\in \Omega^{i-1,j-1}\}.\]

Finally, we denote by $\Omega_p^{i,j}$ the
$\mathrm{SO}(n)$ module of {\it primitive} forms defined as the quotient
\begin{equation} \label{prim1}
\Omega_p^{i,j} :=\Omega^{i,j}/\mathcal{I}^{i,j}.
\end{equation}
An equivalent description of primitive forms can be given as
follows: The multiplication by the symplectic form $-d\alpha$
induces an $\mathrm{SO}(n)$ equivariant linear operator $L:
\Omega_h^{i,j} \rightarrow \Omega_h^{i+1,j+1}$ which is injective
if $i + j \leq n - 2$. Moreover, it follows from the definition
of $\Omega_p^{i,j}$ that in this case
\begin{equation} \label{prim17}
\Omega_p^{i,j} = \Omega_h^{i,j} / L\Omega_h^{i-1,j-1}.
\end{equation}

From Lemma \ref{lemma_induced_reps} and (\ref{prim17}), we now
immediately obtain

\begin{koro} \label{cor_primitive_forms} If $i, j \in \mathbb{N}$
are such that $i+j \leq n-1$, then there is an isomorphism of
$\mathrm{SO}(n)$ modules
\begin{displaymath}
\Omega_p^{i,j} \oplus
\mathrm{Ind}_{\mathrm{SO}(n-1)}^{\mathrm{SO}(n)} (\Lambda^{i-1}
W_0^* \otimes \Lambda^{j-1} W_0^*) \cong
\mathrm{Ind}_{\mathrm{SO}(n-1)}^{\mathrm{SO}(n)}( \Lambda^i W_0^*
\otimes \Lambda^j W_0^*).
\end{displaymath}
\end{koro}

\vspace{0.2cm}

Primitive forms are of particular importance for us since the
space $\mathbf{Val}_i^{\infty}$ fits into an exact sequence of the
spaces $\Omega^{i,j}_p$, as was recently established in
\textbf{\cite{bernig10}}. In order to describe this sequence,
note that $d\mathcal{I}^{i,j} \subseteq \mathcal{I}^{i,j+1}$.
Thus, by definition (\ref{prim1}), the exterior derivative, on
one hand, induces a linear operator $d_Q: \Omega_p^{i,j}
\rightarrow \Omega_p^{i,j+1}$ and, on the other hand, integration
over the normal cycle induces a linear map $\nu:
\Omega^{i,n-i-1}_p \rightarrow \mathbf{Val}_i^{\infty}$ (cf.\
Lemma \ref{smviacurr}). Clearly, both operators are
$\mathrm{SO}(n)$ equivariant.

\begin{theorem} \label{prop_exact_seq}
If $0 \leq i \leq n$, then there is an exact $\mathrm{SO}(n)$
equivariant \linebreak sequence of $\mathrm{SO}(n)$ modules
\begin{displaymath}
0 \to \Lambda^i V_{\mathbb{C}} \rightarrow \Omega^{i,0}_p
\stackrel{\,d_Q}{\rightarrow} \Omega^{i,1}_p
\stackrel{\,d_Q}{\rightarrow} \ldots
\stackrel{\,d_Q}{\rightarrow} \Omega^{i,n-i-1}_p
\stackrel{\nu}{\rightarrow} \mathbf{Val}_i^{\infty} \to 0.
\end{displaymath}
\end{theorem}

\vspace{1cm}

\centerline{\large{\bf{ \setcounter{abschnitt}{5}
\arabic{abschnitt}. Proof of Theorems \ref{thm_decomposition} and
$\mathbf{1'}$}}}

\vspace{0.7cm} \reseteqn \alpheqn \setcounter{theorem}{0}

Theorems \ref{thm_decomposition} and $1'$ are just reformulations
of each other. We give the proof of Theorem
\ref{thm_decomposition} first and then show how Theorem $1'$
follows from it.

\vspace{0.4cm}

\noindent {\it Proof of Theorem \ref{thm_decomposition}}\,: The
cases $i=0$ and $i=n$ are trivial. Moreover, by Theorem
\ref{hardlef}, we may assume that $n/2 \leq i < n$.

Let $\Gamma_{\lambda}$ be an arbitrary irreducible
$\mathrm{SO}(n)$ module of highest weight \linebreak
$\lambda=(\lambda_1,\ldots,\lambda_{\lfloor n/2 \rfloor})$. It is
well known (and a consequence of Corollary
\ref{cor_primitive_forms}) that the multiplicity of
$\Gamma_{\lambda}$ in the spaces $\Omega^{i,j}_p$ of primitive
forms is finite. The same holds true for the spaces
$\mathbf{Val}_i$ since they are quotients of $\Omega^{i,n-i-1}_p$ by Theorem \ref{prop_exact_seq}. Thus, by
Theorem \ref{prop_exact_seq},
we have
\begin{equation} \label{decompval1}
m(\mathbf{Val}_i,\lambda) = (-1)^{n-i}
m(\Lambda^iV_{\mathbb{C}},\lambda) + \sum_{j=0}^{n-i-1}
(-1)^{n-1-i-j} m(\Omega_p^{i,j},\lambda),
\end{equation}
where $m(\,\cdot\,,\lambda)$ denotes the multiplicity of
$\Gamma_{\lambda}$ in the respective $\mathrm{SO}(n)$ modules.

Let $W \cong W^*$ denote the (complex) standard representation of
$\mathrm{SO}(n-1)$. By Corollary \ref{cor_primitive_forms} and
(\ref{propchar1}), we have
\begin{equation*}
m(\Omega^{i,j}_p,\!\lambda)
=m\!\left(\!\mathrm{Ind}_{\mathrm{SO}(n-1)}^{\mathrm{SO}(n)}(\Lambda^iW\!
\otimes \Lambda^jW),\!\lambda\right) - m\!\left(\!
\mathrm{Ind}_{\mathrm{SO}(n-1)}^{\mathrm{SO}(n)}(\Lambda^{i-1}W\!
\otimes \Lambda^{j-1}W),\!\lambda\right)\!.
\end{equation*}
Thus, it follows from an application of Corollary
\ref{cor_clebsch_gordan} (with $n$ replaced by $n-1$ and $0 \leq j
\leq n-1-i$) that
\begin{equation} \label{decompval2}
m(\Omega^{i,j}_p,\lambda) = \sum_{\sigma}
m\!\left(\mathrm{Ind}_{\mathrm{SO}(n-1)}^{\mathrm{SO}(n)}
\bar{\Gamma}_{\sigma},\lambda\right),
\end{equation}
where the sum ranges over all $k:=\lfloor (n-1)/2 \rfloor$-tuples
of non-negative highest weights
$\sigma=(\sigma_1,\ldots,\sigma_k)$ of $\mathrm{SO}(n-1)$ such
that
\[\sigma_1 \leq 2, \qquad \#(\sigma,1)=n-1-i-j, \qquad \#(\sigma,2) \leq j.    \]

If $\mathcal{P}_i$ denotes the union of these $k$-tuples of
non-negative highest weights of $\mathrm{SO}(n-1)$, then, by
(\ref{decompval1}) and (\ref{decompval2}),
\begin{equation} \label{eq_alt_sum}
m(\mathbf{Val}_i,\lambda) = (-1)^{n-i}
m(\Lambda^iV_{\mathbb{C}},\lambda) + \sum_{\sigma \in
\mathcal{P}_i} (-1)^{|\sigma|}
m\!\left(\mathrm{Ind}_{\mathrm{SO}(n-1)}^{\mathrm{SO}(n)}
\bar{\Gamma}_{\sigma},\lambda\right).
\end{equation}

Let $\lambda^*=(\lambda_1^*,\ldots,\lambda_{\lfloor n/2
\rfloor}^*)$, where $\lambda_1^*:=\min\{\lambda_1,2\}$ and
$\lambda_j^*:=|\lambda_j|$ for every $1< j \leq \lfloor n/2
\rfloor$. By Theorem \ref{frobrec}, Theorem \ref{thm_branching}
and the definition of $\bar{\Gamma}_{\sigma}$, we have
\[\sum_{\sigma \in \mathcal{P}_i}
(-1)^{|\sigma|}
m\!\left(\mathrm{Ind}_{\mathrm{SO}(n-1)}^{\mathrm{SO}(n)}
\bar{\Gamma}_{\sigma},\lambda\right)=\sum_{\mu} (-1)^{|\mu|}, \]
where the sum on the right ranges over all sequences
$\mu=(\mu_1,\ldots,\mu_k)$ with $\mu_{n-i}=0$ and
\begin{equation*}
\left \{\begin{array}{ll}  \lambda_1^* \geq \mu_1 \geq \lambda^*_2 \geq \mu_2 \geq \ldots \geq \mu_{k-1} \geq \lambda^*_{\lfloor n/2 \rfloor} \geq |\mu_k| & \quad \mbox{for odd }n, \\
 \lambda_1^* \geq \mu_1 \geq \lambda^*_2 \geq \mu_2 \geq \ldots \geq \mu_{k} \geq \lambda^*_{n/2} & \quad \mbox{for even }n.
\end{array} \right .
\end{equation*}

If $\lambda_{n - i + 1}^*>0$, there is no such sequence. If
$\lambda_{n - i + 1}^*=0$, we obtain
\begin{displaymath}
\sum_{\sigma \in \mathcal{P}_i} (-1)^{|\sigma|}
m\!\left(\mathrm{Ind}_{\mathrm{SO}(n-1)}^{\mathrm{SO}(n)}
\bar{\Gamma}_{\sigma},\lambda\right)=\prod_{j=1}^{n - i - 1}
\sum_{\mu_j=\lambda_{j+1}^*}^{\lambda_j^*} (-1)^{\mu_j}.
\end{displaymath}

This product vanishes if the $\lambda_j^*, j=1,\ldots,n - i$, do
not all have the same parity. If the $\lambda_j^*, j=1,\ldots,n -
i$, all have the same parity, the product above equals $(-1)^{(n
- i-1)\lambda_1^*}$. Consequently, we obtain for $i > n/2$,
\begin{displaymath}
\sum_{\sigma \in \mathcal{P}_i} (-1)^{|\sigma|}
m\!\left(\mathrm{Ind}_{\mathrm{SO}(n-1)}^{\mathrm{SO}(n)}
\bar{\Gamma}_{\sigma},\lambda\right)= \left \{ \begin{array}{ll}
(-1)^{n - i - 1} & \mbox{if } \Gamma_{\lambda} \cong \Lambda^{n - i} V_{\mathbb{C}}, \\
1 & \mbox{if } \lambda \mbox{ satisfies (i), (ii), (iii)},
\\ 0 & \mbox{otherwise.}
\end{array} \right .
\end{displaymath}
If $i=n/2$, in which case $n$ is even, then
\[\sum_{\sigma \in \mathcal{P}_i} (-1)^{|\sigma|}
m\!\left(\mathrm{Ind}_{\mathrm{SO}(n-1)}^{\mathrm{SO}(n)}
\bar{\Gamma}_{\sigma},\lambda\right)= \left \{ \begin{array}{ll}
(-1)^{i-1} & \mbox{if } \lambda = (1,\ldots,1,\pm 1), \\
1 & \mbox{if } \lambda \mbox{ satisfies (i), (ii) and (iii)},
\\ 0 & \mbox{otherwise.}
\end{array} \right .  \]

Plugging this into \eqref{eq_alt_sum} and using that
$\Lambda^{n/2}V_{\mathbb{C}} = \Gamma_{(1,\ldots,1)} \oplus
\Gamma_{(1,\ldots,1,-1)}$ if $n$ is even and
$\Lambda^{n-i}V_{\mathbb{C}} \cong \Lambda^{i}V_{\mathbb{C}}$ for
every $i \in \{0, \ldots, n\}$, completes the proof. \hfill
$\blacksquare$

\pagebreak

Next we explain how Theorem $1'$ can be deduced from Theorem
\ref{thm_decomposition}. The argument presented here in fact shows
that Theorem \ref{thm_decomposition} and $1'$ are equivalent.

\vspace{0.3cm}

\noindent {\it Proof of Theorem $1'$}\,: Let $\Gamma=\Gamma_\mu$
be an irreducible $\mathrm{SO}(n)$-module. The space of $\Gamma$-valued valuations is isomorphic to $\mathbf{Val} \otimes \Gamma$.

Let $S$ denote the
set of highest
weights of $\mathrm{SO}(n)$ satisfying conditions (i)-(iii). By
Theorem \ref{thm_decomposition}, we have
\[\dim (\mathbf{Val}_i \otimes \Gamma)^{\mathrm{SO}(n)} =  \dim (\mathbf{Val}^f_i \otimes \Gamma)^{\mathrm{SO}(n)}=
\sum_{\lambda \in S} \dim(\Gamma_{\lambda} \otimes \Gamma_{\mu})^{\mathrm{SO}(n)}.
\]
Here and in the following, the superscript $\mathrm{SO}(n)$
denotes the subspaces of $\mathrm{SO}(n)$ invariant elements. The
$\Gamma_{\lambda}, \lambda \in S$ are not necessarily self dual
(compare Lemma \ref{lemma_real_reps}). However, if $\lambda \in
S$, then also $\lambda' \in S$, where $\lambda'$ is the highest
weight of $\Gamma_{\lambda}^*$. Thus, by Schur's lemma, we have
\[\dim (\mathbf{Val}_i \otimes \Gamma)^{\mathrm{SO}(n)} = \sum_{\lambda \in S}
\dim \mathrm{Hom}_{\mathrm{SO}(n)}(\Gamma_{\lambda},\Gamma_{\mu})
= \left\{\begin{array}{c l} 1 & \text{if } \mu \in S;\\0 &
\text{otherwise.}\end{array}\right.\]

\vspace{-0.4cm}

\hfill $\blacksquare$

\vspace{0.3cm}

\noindent {\bf Examples:}

\begin{enumerate}
\item[(a)] If $\Gamma=\Gamma_{(0,\ldots,0)} \cong \mathbb{C}$ is the trivial
representation, then $(\mathbf{Val} \otimes \Gamma)^{\mathrm{SO}(n)}\cong\mathbf{Val}^{\mathrm{SO}(n)}$ is the
vector space of all continuous rigid motion invariant valuations
and Theorem $1'$ reduces to Hadwiger's characterization of
intrinsic volumes.

\item[(b)]  If $\Gamma=\Gamma_{(1,0,\ldots,0)}\cong V_{\mathbb{C}}$ is the standard representation, then there is
no
translation invariant and $\mathrm{SO}(n)$ equivariant continuous
valuation with values in $\Gamma$.

\item[(c)] By (\ref{eq_decomp_symmetric}), we have, for $1
\leq i \leq n - 1$,
\[\dim (\mathbf{Val}_i \otimes \mathrm{Sym}^k V_{\mathbb{C}})^{\mathrm{SO}(n)}=\left \{\begin{array}{ll} k/2+1 & \mbox{if $k$ is even} \\
(k-1)/2\quad & \mbox{if $k$ is odd.}   \end{array}  \right .  \]
In particular, if $k = 2$, then there exist (up to constant
multiples) two translation invariant and $\mathrm{SO}(n)$
equivariant continuous $\mathrm{Sym}^2 V_{\mathbb{C}}$ valued
valuations of a given degree $1 \leq i \leq n - 1$. These
valuations are explicitly known (see \textbf{\cite{McMullen97, hug_schneider_schuster_a, Ludwig:matrix}}).

\item[(d)] If $\Gamma=\Gamma_{(1,1,0\ldots,0)}=\Lambda^2 V_{\mathbb{C}}$ is the space of skew-symmetric
tensors of rank two,
then there is no translation invariant and $\mathrm{SO}(n)$
equivariant continuous valuation with values in $\Gamma$. This
answers (the translation invariant case of) a question by Yang
\textbf{\cite{Yang10}}.

\item[(e)] The unique translation invariant and $\mathrm{SO}(n)$
equivariant continuous \linebreak valuation with values in
$\Gamma_{(2,\ldots,2,0,\ldots,0)}$ was constructed in
\textbf{\cite{Bernig06}}.
\end{enumerate}

\vspace{0.8cm}

\centerline{\large{\bf{ \setcounter{abschnitt}{6}
\arabic{abschnitt}. Bivaluations}}}

\vspace{0.6cm} \reseteqn \alpheqn \setcounter{theorem}{0}

We turn now to the study of bivaluations. In particular, we will
present the proof of Theorem \ref{thm_symmetry} at the end of this section.

We denote the vector space of all continuous translation
biinvariant \linebreak complex valued bivaluations by
$\mathbf{BVal}$ and we write $\mathbf{BVal}_{i,j}$ for its
subspace of all bivaluations of bidegree $(i,j)$. An immediate
consequence of McMullen's decomposition \eqref{eq_mcmullen} of
the vector space $\mathbf{Val}$ is a corresponding result for the
space $\mathbf{BVal}$:

\begin{koro} \label{mcmullenbval} The space $\mathbf{BVal}$ admits a decomposition
\[\mathbf{BVal}=\bigoplus_{i,j=0}^n \mathbf{BVal}_{i,j}.\]
\end{koro}

Corollary \ref{mcmullenbval} implies an analog of
Corollary \ref{normval} for the space of translation biinvariant
bivaluations as follows.

\begin{koro} Let $C \in \mathcal{K}^n$ be a fixed convex body with
non-empty interior. The space $\mathbf{BVal}$ becomes a Banach
space under the norm
\[\|\varphi \| = \sup \{|\varphi(K,L)|: K, L \subseteq C \}.  \]
Moreover, a different choice of $C$ yields an equivalent norm.
\end{koro}

The group $\mathrm{O}(n)\times \mathrm{O}(n)$ acts continuously
on the Banach space $\mathbf{BVal}$ by
\[((\eta,\vartheta)\phi)(K,L)=\phi(\eta^{-1} K,\vartheta^{-1}L), \qquad
(\eta,\vartheta) \in \mathrm{O}(n)\times \mathrm{O}(n), \varphi
\in \mathbf{BVal}.\]

We denote by $\mathbf{BVal}^f$ the subspace of bivaluations with
finite dimensional $\mathrm{O}(n)\times \mathrm{O}(n)$ orbit.
Since $\mathrm{O}(n)\times \mathrm{O}(n)$ is compact,
$\mathbf{BVal}^f$ is a dense subspace \linebreak of
$\mathbf{BVal}$ (see e.g.\, \textbf{\cite[\textnormal{p.\
141}]{broecker_tomdieck}}).

\begin{prop} \label{prop_isom_bival}
Let $0 \leq i,j \leq n$. The linear map $\iota: \mathbf{Val}_i^f
\otimes \mathbf{Val}_j^f \to \mathbf{BVal}_{i,j}^f$, induced by
\begin{equation} \label{eq_map_bivaluations}
\iota(\phi \otimes \psi)(K,L)=\phi(K)\psi(L),
\end{equation}
is an isomorphism of  $\mathrm{O}(n)\times \mathrm{O}(n)$ modules.
\end{prop}

\noindent {\it Proof}\,: It is easy to see that the map $\iota$ is
$\mathrm{O}(n)\times \mathrm{O}(n)$ equivariant and injective. It
remains to prove that it is onto.

It is well known that every irreducible $\mathrm{O}(n) \times
\mathrm{O}(n)$ module is of the form $\Gamma \otimes \Theta$,
where $\Gamma, \Theta$ are irreducible $\mathrm{O}(n)$ modules
(c.f.\ \textbf{\cite[\textnormal{p.\ 82}]{broecker_tomdieck}}).
Thus, if \linebreak $\varphi \in \mathbf{BVal}^f_{i,j}$ belongs to
a subspace isomorphic to $\Gamma \otimes \Theta$, the valuation
$\varphi(\,\cdot\,,L) \in \mathbf{Val}_i$ belongs to a subspace
which is isomorphic to $\Gamma$ as an $\mathrm{O}(n)$ module for
every $L \in \mathcal{K}^n$. Since any $\mathrm{O}(n)$
representation whose restriction \linebreak to $\mathrm{SO}(n)$ is
multiplicity free, is itself multiplicity free, it follows from
\linebreak Theorem \ref{thm_decomposition} that the irreducible
subspace of $\mathbf{Val}_i$ which is isomorphic to $\Gamma$ has
multiplicity at most one.

If $\{\phi_1,\ldots,\phi_l\}$ is a basis of the isomorphic copy of
$\Gamma$ in $\mathbf{Val}_i$, then
\[\varphi(\,\cdot\,,L)=\sum_{k=1}^l\phi_k(\,\cdot\,)\psi_k(L),\]
where $\psi_k(L)$ are coefficients depending on $L$. It is not
difficult to show that $\psi_k \in \mathbf{Val}_j$ and that
$\psi_k$ belongs to an isomorphic copy of $\Theta$ in
$\mathbf{Val}_j$ for every $k \in \{1, \ldots, l\}$.
Thus, we have
shown that $\varphi$ is the image under the map $\iota$ of the
element $\sum_{k=1}^l \phi_k\otimes \psi_k \in \mathbf{Val}_i^f \otimes \mathbf{Val}_j^f$. \hfill
$\blacksquare$

\vspace{0.5cm}

After these preparations, we are now in a position to proof the
following refinement of Theorem \ref{thm_symmetry}.

\begin{theorem} Suppose that $\varphi \in \mathbf{BVal}_{i,i}$, where
$0 \leq i \leq n$.
\begin{enumerate}
\item[(a)] If $\varphi$ is $\mathrm{O}(n)$ invariant, then
$\varphi(K,L)=\varphi(L,K)$ for every $K, L \in \mathcal{K}^n$.
\item[(b)] If $\varphi$ is $\mathrm{SO}(n)$ invariant and $(i,n)
\neq (2k + 1,4k + 2), k \in \mathbb{N}$, then
$\varphi(K,L)=\varphi(L,K)$ for every $K, L \in \mathcal{K}^n$.
\item[(c)] If $(i,n) = (2k + 1,4k + 2), k \in \mathbb{N}$, then
there exist an $\mathrm{SO}(n)$ invariant $\zeta \in
\mathbf{BVal}_{i,i}$ and $K, L \in \mathcal{K}^n$ such that
$\zeta(K,L) \neq \zeta(L,K)$.

\end{enumerate}
\end{theorem}

\noindent {\it Proof}\,: Since the cases $i = 0$ and $i = n$ are
trivial, we may assume that $0<i<n$. Moreover, since
$\mathrm{O}(n) \times \mathrm{O}(n)$ finite bivaluations are
dense in $\mathbf{BVal}_{i,i}$ we may assume that $\varphi \in
\mathbf{BVal}_{i,i}^f$, where $1 \leq i \leq n$.

From Theorem \ref{thm_decomposition} we deduce that the
decom\-position of the space $\mathbf{Val}_i$ into irreducible
$\mathrm{O}(n)$ modules is multiplicity free, say
\begin{displaymath}
\mathbf{Val}_i=\bigoplus_{\gamma \in R} \Gamma_\gamma,
\end{displaymath}
where the sum ranges over some set $R$ of equivalence classes of
irreducible representations of $\mathrm{O}(n)$.

From Proposition \ref{prop_isom_bival}, it follows that
\[\mathbf{BVal}_{i,i}^{f,\mathrm{O}(n)} \cong (\mathbf{Val}_i^f \otimes \mathbf{Val}_i^f)^{\mathrm{O}(n)}
 \cong \bigoplus_{\gamma,\delta \in R} (\Gamma_{\gamma} \otimes
\Gamma_{\delta})^{\mathrm{O}(n)}.\] Since, by Lemma
\ref{lemma_real_reps}, all representations of $\mathrm{O}(n)$ are
self-dual, we have
\[(\Gamma_{\gamma} \otimes \Gamma_{\delta})^{\mathrm{O}(n)} \cong \mathrm{Hom}_{\mathrm{O}(n)}(\Gamma_{\gamma},\Gamma_{\delta})
 \cong \mathrm{Hom}_{\mathrm{O}(n)}(\Gamma_{\gamma} \otimes \Gamma_{\delta},\mathbb{C}).  \]
Since $\Gamma_{\gamma}$ and $\Gamma_{\delta}$ are irreducible,
Schur's lemma implies that
\[\dim \mathrm{Hom}_{\mathrm{O}(n)}(\Gamma_{\gamma},\Gamma_{\delta})=
\left\{ \begin{array}{l l} 1 \quad & \mbox{if } \gamma = \delta
\\0 & \mbox{if } \gamma \neq \delta. \end{array}\right.\]
Since all representations of $\mathrm{O}(n)$ are real, the space
\[\mathrm{Hom}_{\mathrm{O}(n)}(\Gamma_{\gamma} \otimes
\Gamma_{\gamma},\mathbb{C})=(\mathrm{Sym}^2\Gamma_{\gamma})^{\mathrm{O}(n)}
\oplus (\mathrm{\Lambda}^2\Gamma_{\gamma})^{\mathrm{O}(n)}
\] of $\mathrm{O}(n)$ invariant bilinear forms on
$\Gamma_{\gamma}$ coincides with
$(\mathrm{Sym}^2\Gamma_{\gamma})^{\mathrm{O}(n)}$. Thus,
\[\mathbf{BVal}_{i,i}^{f,\mathrm{O}(n)} \cong \bigoplus_{\gamma \in R}\, (\mathrm{Sym}^2\Gamma_{\gamma} )^{\mathrm{O}(n)}\]
which completes the proof of statement (a).

\pagebreak

If $n \not \equiv 2 \mod 4$, then the proof of statement (b) is
similar, since in this case all representations of
$\mathrm{SO}(n)$ are also real, by Lemma \ref{lemma_real_reps}.
However, in the case $n \equiv 2 \mod 4$, more care is needed,
since there are $\mathrm{SO}(n)$ modules which are not real. By
Lemma \ref{lemma_real_reps}, an irreducible $\mathrm{SO}(n)$
module $\Gamma_\lambda$ of highest weight $\lambda =
(\lambda_1,\ldots,\lambda_{n/2})$ is real if and only if
$\lambda_{n/2} = 0$. If $i \neq n/2$, then, by Theorem
\ref{thm_decomposition}, all irreducible $\mathrm{SO}(n)$ modules
which enter $\mathbf{Val}_i$ are of this form. Consequently, any
$\mathrm{SO}(n)$ invariant bivaluation $\varphi \in
\mathbf{BVal}_{i,i}$ is symmetric in this case.

Finally let $n \equiv 2 \mod 4$ and $i = n/2$. By Theorem
\ref{thm_decomposition}, the dual irreducible $\mathrm{SO}(n)$
modules $\Gamma_{(2,\ldots,2)}$ and $\Gamma_{(2,\ldots,2,-2)}$
both enter $\mathbf{Val}_i$ with multiplicity one. If
$\{\phi_1,\ldots,\phi_l\}$ is a basis of $\Gamma_{(2,\ldots,2)}
\subseteq \mathbf{Val}_i^f$ and $\{\psi_1,\ldots,\psi_l\}$ denotes
the corresponding dual basis in $\Gamma_{(2,\ldots,2,-2)}
\subseteq \mathbf{Val}_i^f$, then the image of $\sum_{k=1}^l
\phi_k \otimes \psi_k$ under the map $\iota$ defined in
\eqref{eq_map_bivaluations} clearly is a continuous
$\mathrm{SO}(n)$ invariant bivaluation in $\mathbf{BVal}_{i,i}$.
However, it is not symmetric since the valuations
$\{\phi_1,\ldots,\phi_l, \psi_1,\ldots,\psi_l\}$ are linearly
independent. \hfill $\blacksquare$

\vspace{1cm}

\centerline{\large{\bf{ \setcounter{abschnitt}{7}
\arabic{abschnitt}. Applications to geometric inequalities}}}

\vspace{0.7cm} \reseteqn \alpheqn \setcounter{theorem}{0}

As applications of Theorem \ref{thm_symmetry}, we present in this section several
new \linebreak geometric inequalities involving $\mathrm{SO}(n)$
equivariant Minkowski valuations. Their proofs are based, on one
hand, on the symmetry of bivaluations and, on the other hand, on
techniques developed by Lutwak \textbf{\cite{lutwak85, lutwak86,
lutwak88, lutwak90, lutwak93}}.

\begin{lem} \label{sonon} If $\Phi \in \mathbf{MVal}$ is $\mathrm{SO}(n)$
equivariant, then $\Phi$ is also $\mathrm{O}(n)$ \linebreak
equivariant.
\end{lem}

\noindent {\it Proof}\,: Let $\mathbf{CVal}$ denote the vector
space of all continuous and translation invariant valuations with
values in the space $C(S^{n-1})$ of continuous complex valued
functions on $S^{n-1}$. Note that any $\mathrm{SO}(n)$
equivariant $\Phi \in \mathbf{MVal}$ induces an $\mathrm{SO}(n)$
equivariant $\bar{\Phi} \in \mathbf{CVal}$, by
$\bar{\Phi}(K,\cdot) = h(\Phi K,\cdot)$. Therefore, it is
sufficient to show that any $\mathrm{SO}(n)$ equivariant
valuation in $\mathbf{CVal}$ is $\mathrm{O}(n)$ equivariant.

Using arguments as in the proof of Proposition
\ref{prop_isom_bival}, it is easy to show that $\mathbf{CVal}^{f}
\cong \mathbf{Val}^f \otimes C(S^{n-1})^f$ as $\mathrm{O}(n)
\times \mathrm{O}(n)$ modules. Consequently,
\begin{equation} \label{cval}
\mathbf{CVal}^{f,\mathrm{O}(n)} \cong (\mathbf{Val}^f \otimes
C(S^{n-1})^f)^{\mathrm{O}(n)}.
\end{equation}

It is well known that the decomposition of $C(S^{n-1})$ into
irreducible $\mathrm{SO}(n)$ modules is given by
\begin{equation} \label{csdecomp}
C(S^{n-1})= \bigoplus_{k \geq 0} \Gamma_{(k,0,\ldots,0)},
\end{equation}
where the spaces $\Gamma_{(k,0,\ldots,0)}$ are precisely the
spaces of spherical harmonics of degree $k$ in dimension $n$.
Moreover, the spaces $\Gamma_{(k,0,\ldots,0)}$ are self-dual and
$\mathrm{O}(n)$ invariant and, thus, (\ref{csdecomp}) also
represents the decomposition of $C(S^{n-1})$ into irreducible
$\mathrm{O}(n)$ modules.

Let $m_k$ denote the (finite) multiplicity of the isomorphic copy
of $\Gamma_{(k,0,\ldots,0)}$ in $\mathbf{Val}$. From
(\ref{cval}), Theorem \ref{thm_decomposition} and an application
of Schur's lemma, we obtain
\[\mathbf{CVal}^{\mathrm{SO}(n)}=
\bigoplus_{k} m_k (\Gamma_{(k,0,\ldots,0)} \otimes
\Gamma_{(k,0,\ldots,0)})^{\mathrm{SO}(n)} =
\mathbf{CVal}^{\mathrm{O}(n)}.  \] Thus, any $\mathrm{SO}(n)$
equivariant valuation in $\mathbf{CVal}$ is also $\mathrm{O}(n)$
equivariant. \hfill $\blacksquare$

\vspace{0.4cm}

For $K, L \in \mathcal{K}^n$ and $0 \leq i \leq n - 1$, we write
$W_i(K,L)$ to denote the mixed volume
$V(K,\ldots,K,B,\ldots,B,L)$, where $K$ appears $n-1-i$ times and
the Euclidean unit ball $B$ appears $i$ times. The mixed volume
$W_i(K,K)$ will be written as $W_i(K)$ and is called the
\emph{$i$-th quermassintegral} of $K$. The \emph{$i$-th intrinsic
volume} $V_i(K)$ of $K$ is defined by
\begin{equation} \label{viwi}
\kappa_{n-i}V_i(K)=\binom{n}{i} W_{n-i}(K),
\end{equation}
where $\kappa_n$ is the $n$-dimensional volume of the Euclidean
unit ball in $V$.

We will repeatedly make use of the following consequence
of Theorem \ref{thm_symmetry} and Lemma \ref{sonon}.

\begin{koro} \label{symkor} If $\Phi_i \in \mathbf{MVal}_i$, $1 \leq i \leq n - 1$, is $\mathrm{SO}(n)$
equivariant, then
\[W_{n-1-i}(K,\Phi_iL) = W_{n-1-i}(L,\Phi_iK)   \]
for every $K, L \in \mathcal{K}^n$. \hfill $\blacksquare$
\end{koro}

\vspace{0.4cm}

Let $\mathcal{K}^n_{\mathrm{o}} \subseteq \mathcal{K}^n$ denote
the set of convex bodies with non-empty interior. \linebreak One
of the fundamental inequalities for mixed volumes is the
(general) Minkowski inequality: If $K, L \in
\mathcal{K}^n_{\mathrm{o}}$ and $0 \leq i \leq n - 2$, then
\begin{equation} \label{genmink}
W_i(K,L)^{n-i} \geq W_i(K)^{n-i-1}W_i(L),
\end{equation}
with equality if and only if $K$ and $L$ are homothetic.

\begin{lem} \label{upbound} Let $\Phi_i \in \mathbf{MVal}_i$, $1 \leq i \leq n - 1$, be $\mathrm{SO}(n)$
equivariant and non-trivial, i.e.\ $\Phi_i(K) \neq \{0\}$ for some
$K \in \mathcal{K}^n$.
\begin{enumerate}
\item[(a)] There exists a constant $r(\Phi_i) > 0$ such that
for every $K \in \mathcal{K}^n$,
\[W_{n-1}(\Phi_iK)=r(\Phi_i)W_{n-i}(K).   \]
\item[(b)] If $K \in \mathcal{K}^n_{\mathrm{o}}$, then
\[W_{n-i}(K)^{i+1} \geq \frac{\kappa_{n}^{i}}{r(\Phi_i)^{i+1}}W_{n-1-i}(\Phi_i K).   \]
If $\Phi_i \mathcal{K}^n_{\mathrm{o}} \subseteq
\mathcal{K}^n_{\mathrm{o}}$, then equality holds if and only if
$\Phi_iK$ is a ball.
\end{enumerate}
\end{lem}

\noindent {\it Proof}\,: Statement (a) follows from Hadwiger's
characterization theorem. From repeated application of
Minkowski's inequality (\ref{genmink}) with $L=B$, we obtain the
inequality
\[W_{n-1}(K)^{i+1} \geq \kappa_n^iW_{n-1-i}(K),   \]
where, for $K \in \mathcal{K}^n_{\mathrm{o}}$, there is equality
if and only if $K$ is a ball. Taking $\Phi_i K$ instead of $K$ and using (a), yields statement (b). \hfill $\blacksquare$ \\

Special cases of Lemma \ref{upbound} were previously obtained by
Lutwak \textbf{\cite{lutwak85}} (for $\Phi_i=\Pi_i$) and one of
the authors \textbf{\cite{Schu06a}} (for $i=n-1$).

In order to proof Theorem \ref{thm_brunn_minkowski}, we need a
further generalization of the Brunn--Minkowski inequality
\eqref{intvolbm} (where the equality conditions are not yet
known): If $0 \leq i \leq n-2$, $K, L, K_1, \ldots, K_i \in
\mathcal{K}^n$ and $\mathbf{C}=(K_1,...,K_i)$, then
\begin{equation} \label{mostgenbm}
V_i(K+L,\mathbf{C})^{1/(n-i)} \geq
V_i(K,\mathbf{C})^{1/(n-i)}+V_i(L,\mathbf{C})^{1/(n-i)}.
\end{equation}

{\it Proof of Theorem \ref{thm_brunn_minkowski}}\,: Since
translation invariant continuous Minkowski valuations which are
homogeneous of degree one are linear with respect to Minkowski
addition (see e.g.\ \textbf{\cite{hadwiger51}}), the case $i=1$
is a direct consequence of inequality (\ref{intvolbm}). Thus, we
may assume that $i \geq 2$.

By Corollary \ref{symkor} and (\ref{mostgenbm}), we
have for $Q \in \mathcal{K}^n_{\mathrm{o}}$,
\begin{eqnarray*}
W_{n-1-i}(Q,\Phi_i(K+L))^{1/i} & = & W_{n-1-i}(K+L,\Phi_iQ)^{1/i}
\\
& \geq & W_{n-1-i}(K,\Phi_i Q)^{1/i}+W_{n-1-i}(L,\Phi_i Q)^{1/i}
\\
&=& W_{n-1-i}(Q,\Phi_i K)^{1/i}+W_{n-1-i}(Q,\Phi_i L)^{1/i}.
\end{eqnarray*}
It follows from Minkowski's inequality (\ref{genmink}), that
\begin{equation} \label{mink1}
W_{n-1-i}(Q,\Phi_iK)^{i+1}\geq W_{n-1-i}(Q)^{i}W_{n-1-i}(\Phi_iK),
\end{equation}
and
\begin{equation} \label{mink2}
W_{n-1-i}(Q,\Phi_iL)^{i+1}\geq W_{n-1-i}(Q)^{i}W_{n-1-i}(\Phi_iL).
\end{equation}
Thus, if we set $Q=\Phi_i(K+L)$ and use \eqref{viwi}, we obtain
the desired inequality
\begin{equation*}
W_{n-1-i}(\Phi_i(K+L))^{1/i(i+1)} \geq
W_{n-1-i}(\Phi_iK)^{1/i(i+1)}+W_{n-1-i}(\Phi_iL)^{1/i(i+1)}.
\end{equation*}

Suppose now that equality holds and that $\Phi_i
\mathcal{K}^n_{\mathrm{o}} \subseteq \mathcal{K}^n_{\mathrm{o}}$.
By Theorem $1'$, applied to the standard representation
$V=\Gamma_{(1,0,\ldots,0)}$, the Steiner point of $\Phi_iK$ is
the origin for every $K \in \mathcal{K}^n$. Thus, we can deduce
from the equality conditions of (\ref{mink1}) and (\ref{mink2}),
that there exist $\lambda_1, \lambda_2 > 0$ such that
\begin{equation} \label{hom2}
\Phi_iK=\lambda_1\Phi_i(K+L) \qquad \mbox{and} \qquad
\Phi_iL=\lambda_2 \Phi_i(K+L)
\end{equation}
and
\[\lambda_1^{1/i}+\lambda_2^{1/i} = 1.   \]
Using Lemma \ref{upbound} (a) and (\ref{hom2}) we get
\[W_{n-i}(K+L)^{1/i}=W_{n-i}(K)^{1/i}+W_{n-i}(L)^{1/i},  \]
which, by (\ref{intvolbm}), implies that $K$ and $L$ are
homothetic.
\hfill $\blacksquare$\\

The major open problem concerning the rigid motion invariant
quantities $W_{n-1-i}(\Phi_iK)$ is how to estimate them from
below in terms of $W_{n-1-i}(K)$. A standard method of proof for
isoperimetric problems of this kind was introduced by Lutwak
\textbf{\cite{lutwak86}} and is now known as the {\it class
reduction technique}. Our last result shows how Corollary
\ref{symkor} allows for applications of the class reduction
technique to the functionals $W_{n-1-i}(\Phi_iK)$, $K \in
\mathcal{K}^n$.

In the following we use $\Phi_i^2K$ to denote $\Phi_i\Phi_iK$.

\begin{theorem} \label{classred} Let $\Phi_i \in \mathbf{MVal}_i$, $1 \leq i \leq n - 1$, be $\mathrm{SO}(n)$
equivariant and suppose that $\Phi_i \mathcal{K}^n_{\mathrm{o}}
\subseteq \mathcal{K}^n_{\mathrm{o}}$. If $K \in
\mathcal{K}^n_{\mathrm{o}}$, then
\[\frac{W_{n-1-i}(\Phi_iK)}{W_{n-1-i}(K)^i} \geq \frac{W_{n-1-i}(\Phi_i^2K)}{W_{n-1-i}(\Phi_iK)^i},   \]
with equality if and only if $K$ and $\Phi_i^2 K$ are homothetic.
\end{theorem}
\noindent {\it Proof}\,: Let $K, Q \in
\mathcal{K}^n_{\mathrm{o}}$. From Corollary \ref{symkor} and the
Minkowski inequality (\ref{genmink}), we obtain
\[W_{n-1-i}(Q,\Phi_i K)^{i+1}=W_{n-1-i}(K,\Phi_i Q)^{i+1}\geq W_{n-1-i}(K)^iW_{n-1-i}(\Phi_i Q),    \]
with equality if and only if $K$ and $\Phi_i Q$ are homothetic.
Taking $Q=\Phi_i K$, yields
\[W_{n-1-i}(\Phi_i K)^{i+1} \geq W_{n-1-i}(K)^{i}W_{n-1-i}(\Phi_i^2K),  \]
with equality if and only if $K$ and $\Phi_i^2 K$ are homothetic.
\hfill $\blacksquare$ \\

\begin{small}

{\scshape Tel Aviv University, Israel} \par {\it E-mail address}: semyon@post.tau.ac.il \\

{\scshape Goethe-Universit\"at Frankfurt am Main, Germany} \par {\it E-mail address}:
bernig@math.uni-frankfurt.de  \\

{\scshape Vienna University of Technology, Austria} \par {\it
E-mail address}: franz.schuster@tuwien.ac.at

\end{small}


\vspace{0.1cm}




\begin{thebibliography}{99}
\footnotesize{
\parskip-0.1cm{

\bibitem{Abardia_Bernig11} J. Abardia, A. Bernig, \emph{Projection bodies in complex vector spaces}, Adv. Math., in press.

\bibitem{adams} J. F. Adams, \emph{Lectures on {L}ie groups}, W. A. Benjamin, Inc., 
New York-Amsterdam, 1969.

\bibitem{Alesker99}
S. Alesker, \emph{Continuous rotation invariant valuations on
convex sets}, Ann. of Math. (2) {\bf 149} (1999), 977--1005.

\bibitem{Alesker97} S. Alesker, \emph{Description of continuous isometry covariant
valuations on convex sets}, Geom. Dedicata {\bf 74} (1999),
241--248.

\bibitem{Alesker01}
S. Alesker, \emph{Description of translation invariant valuations
on convex sets with solution of P. McMullen's conjecture}, Geom.
Funct. Anal. {\bf 11} (2001), 244--272.

\bibitem{Alesker03}
S. Alesker, \emph{Hard Lefschetz theorem for valuations, complex
integral geometry, and unitarily invariant valuations}, J.
Differential Geom. {\bf 63} (2003), 63--95.

\bibitem{Alesker04b}
S. Alesker, \emph{The multiplicative structure on polynomial
valuations}, Geom. Funct. Anal. {\bf 14} (2004), 1--26.

\bibitem{Alesker06a}
S. Alesker, \emph{Theory of valuations on manifolds. I. Linear
spaces}, Israel J. Math.  {\bf 156} (2006), 311--339.

\bibitem{alesker10}
S. Alesker, \emph{A Fourier type transform on translation
invariant valuations on convex sets}, Israel J. Math., in press.

\bibitem{AlBern}
S. Alesker and J. Bernstein, \emph{Range characterization of the
cosine transform on higher Grassmannians}, Adv. Math. {\bf 184}
(2004), 367--379.

\bibitem{Bernig06}
A. Bernig, \emph{Curvature tensors of singular spaces}, Diff. Geom. Appl. {\bf 24} (2006), 191--208.

\bibitem{Bernig09}
A. Bernig, \emph{A Hadwiger-type theorem for the special unitary
group}, Geom. Funct. Anal. {\bf 19} (2009), 356--372.

\bibitem{bernig10} A. Bernig, \emph{Invariant valuations on quaternionic vector spaces},
J. Inst. de Math. de Jussieu, in press.

\bibitem{Bernig07b}
A. Bernig and L. Br\"ocker, \emph{Valuations on manifolds and
Rumin cohomology},  J. Differential Geom.  {\bf 75} (2007),
433--457.

\bibitem{bernigfu06}
A. Bernig and J.H.G. Fu, \emph{Convolution of convex valuations},
Geom. Dedicata {\bf 123} (2006), 153--169.

\bibitem{bernigfu10}
A. Bernig and J.H.G. Fu, \emph{Hermitian integral geometry}, Ann.
of Math. (2), in press.

\bibitem{broecker_tomdieck}
T. Br\"ocker, T. tom Dieck, \emph{Representations of compact Lie
groups}. Graduate Texts in Mathematics {\bf 98}, Springer-Verlag,
New York, 1985.

\bibitem{fu06}
J.H.G. Fu, \emph{Structure of the unitary valuation algebra}, J.
Differential Geom. {\bf 72} (2006), 509--533.

\bibitem{fultonharris}
W. Fulton and J. Harris, \emph{Representation theory. A first
course}, Graduate Texts in Mathematics {\bf 129}, Springer-Verlag,
New York, 1991.

\bibitem{goodeyzhang}
P. Goodey and G. Zhang, \emph{Inequalities between projection
functions of convex bodies}, Amer. J. Math. {\bf 120} (1998),
345--367.

\bibitem{goodwallach}
R. Goodman and N.R. Wallach, \emph{Representations and invariants
of the classical groups}, Encyclopedia of Mathematics and its
Applications {\bf 68}, Cambridge University Press, Cambridge,
1998.

\bibitem{grinberg:zhang}
E. Grinberg and G. Zhang, \emph{Convolutions, transforms, and
convex bodies}, Proc. London Math. Soc. (3) \textbf{78} (1999)
77--115.

\bibitem{haberl08}
C. Haberl, \emph{$L_p$ intersection bodies}, Adv. Math. {\bf 217}
(2008), 2599--2624.

\bibitem{haberl09}
C. Haberl, \emph{Blaschke valuations}, Amer. J. Math., in press.

\bibitem{habschu09}
C. Haberl and F.E. Schuster, \emph{General $L_p$ affine
isoperimetric inequalities}, J. Differential Geom. {\bf 83} (2009), 1--26.

\bibitem{hadwiger51a}
H. Hadwiger, \emph{Translationsinvariante, additive und stetige
Eibereichfunktionale}, Publ. Math. Debrecen {\bf 2} (1951), 81-94.

\bibitem{hadwiger51}
H. Hadwiger, \emph{Vorlesungen \"uber Inhalt, Oberfl\"ache und
Isoperimetrie}, Springer, Berlin, 1957.

\bibitem{hug_schneider_schuster_a}
D. Hug, R. Schneider and R. Schuster, \emph{The space of isometry covariant tensor valuations}, Algebra i Analiz {\bf 19} (2007), 194--224.

\bibitem{hug_schneider_schuster_b}
D. Hug, R. Schneider and R. Schuster, \emph{Integral geometry of tensor valuations}, Adv. in Appl. Math. {\bf 41} (2008), 482--509.

\bibitem{kiderlen05}
M. Kiderlen, \emph{Blaschke- and Minkowski-Endomorphisms of convex
bodies}, Trans. Amer. Math. Soc. {\bf 358} (2006), 5539--5564.

\bibitem{Klain:Rota}
D.A. Klain and G.-C. Rota, \emph{Introduction to geometric
probability}, Cambridge University Press, Cambridge, 1997.

\bibitem{ludwig02} M. Ludwig, \emph{Projection bodies and
valuations}, Adv. Math. {\bf 172} (2002), 158-168.

\bibitem{Ludwig:matrix}
M. Ludwig, \emph{Ellipsoids and matrix valued valuations}, Duke
Math. J. {\bf 119} (2003), 159--188.

\bibitem{Ludwig:Minkowski}
M. Ludwig, \emph{Minkowski valuations}, Trans. Amer. Math. Soc.
\textbf{357} (2005) 4191--4213.

\bibitem{Ludwig06}
M. Ludwig, \emph{Intersection bodies and valuations}, Amer. J.
Math. {\bf 128} (2006), 1409--1428.

\bibitem{Ludwig10a}
M. Ludwig, \emph{Minkowski areas and valuations}, J. Differential
Geom., in press.

\bibitem{centro}
M. Ludwig and M. Reitzner, \emph{A classification of
$\mathrm{SL}(n)$ invariant valuations}, Ann. of Math. (2) {\bf 172} (2010), 1223--1271.

\bibitem{lutwak85}
E. Lutwak, \emph{Mixed projection inequalities}, Trans. Amer.
Math. Soc. {\bf 287} (1985), 91--105.

\bibitem{lutwak86}
E. Lutwak, \emph{On some affine isoperimetric inequalities}, J.
Differential Geom. \textbf{23} (1986), 1--13.

\bibitem{lutwak88} E. Lutwak, \emph{Intersection bodies and dual
mixed volumes}, Adv. Math. {\bf 71} (1988), 232--261.

\bibitem{lutwak90} E. Lutwak, {\em Centroid bodies and dual mixed
volumes}, Proc. London Math. Soc. {\bf 60} (1990), 365--391.

\bibitem{lutwak93}
E. Lutwak, \emph{Inequalities for mixed projection bodies}, Trans.
Amer. Math. Soc. {\bf 339} (1993), no. 2, 901-916.

\bibitem{LYZ2000a}
E. Lutwak, D. Yang, and G. Zhang, \emph{$L_p$ affine
isoperimetric inequalities}, J. Differential Geom. {\bf 56}
(2000), 111--132.

\bibitem{LYZ2000b}
E. Lutwak, D. Yang, and G. Zhang, \emph{A new ellipsoid associated
with convex bodies}, Duke Math. J. {\bf 104} (2000) 375--390.

\bibitem{LZ1997}
E. Lutwak and G. Zhang, \emph{Blaschke-Santal\'o inequalities},
J. Differential Geom. \textbf{47} (1997), 1--16.

\bibitem{McMullen77}
P. McMullen, \emph{Valuations and Euler-type relations on certain
classes of convex polytopes}, Proc. London Math. Soc. {\bf 35}
(1977), 113--135.

\bibitem{McMullen93}
P. McMullen, \emph{Valuations and dissections}, Handbook of Convex
Geometry, Vol.\ B (P.M. Gruber and J.M. Wills, eds.),
North-Holland, Amsterdam, 1993, pp. 933--990.

\bibitem{McMullen97}
P. McMullen, \emph{Isometry covariant valuations on convex
bodies}, Rend. Circ. Mat. Palermo (2) Suppl. {\bf 50} (1997),
259--271.

\bibitem{rumin}
M. Rumin, \emph{Formes diff\'erentielles sur les vari\'et\'es de
contact}, J. Differential Geom. {\bf 39} (1994), 281--330.

\bibitem{Schu06a}
F.E. Schuster, \emph{Convolutions and multiplier
transformations}, Trans. Amer. Math. Soc. {\bf 359} (2007),
5567--5591.

\bibitem{Schu09}
F.E. Schuster, \emph{Crofton Measures and Minkowski Valuations},
Duke Math. J. {\bf 154} (2010), 1--30.

\bibitem{wallach1}
N.R. Wallach, \emph{Real reductive groups. I}, Pure and Applied
Mathematics {\bf 132}, Academic Press, Inc., Boston, MA, 1988.

\bibitem{Yang10} D. Yang, \emph{Affine integral geometry from a differentiable viewpoint},
Handbook of Geometri\-c Analysis, No. 2, in press.

}}
\end{thebibliography}
\end{document}